\documentclass{dcds-bOF}
\usepackage{amsmath}
  \usepackage{paralist}
  \usepackage{graphics} 
  \usepackage{epsfig} 
\usepackage{graphicx}  \usepackage{epstopdf}
 \usepackage[colorlinks=true]{hyperref}
\hypersetup{urlcolor=blue, citecolor=red}

\usepackage{dsfont}

\def\currentvolume{X}
 \def\currentissue{X}
  \def\currentyear{200X}
   \def\currentmonth{XX}
    \def\ppages{X--XX}
     \def\DOI{10.3934/dcdsb.2019154}

\newtheorem{theorem}{Theorem}[section]
\newtheorem{assumption}[theorem]{Assumption}
\newtheorem{corollary}{Corollary}

\newtheorem{lemma}[theorem]{Lemma}

\theoremstyle{definition}

\newtheorem{remark}{Remark}

\newcommand{\N}{\mathbb{N}}
\newcommand{\R}{\mathbb{R}}
\newcommand{\E}{\mathbb{E}}
\renewcommand{\P}{\mathbb{P}}
\newcommand{\1}{\mathds{1}}
\newcommand{\F}{\mathcal{F}}
\newcommand{\tr}{\operatorname{trace}}

\title[Mean-square approximations of L\'{e}vy noise driven SDEs] 
      {Mean-square approximations of L\'{e}vy noise driven SDEs with super-linearly growing diffusion and jump coefficients}

\author[Ziheng Chen, Siqing Gan and Xiaojie Wang]{}


 \email{zihengchen@csu.edu.cn}
 \email{sqgan@csu.edu.cn}
 \email{x.j.wang7@csu.edu.cn}

\thanks{This work was supported by NSF of China (11571373, 11671405, 91630312), NSF of Hunan Province (2016JJ3137), Innovation-Driven Project of CSU (2017CX017),  Shenghua Yuying Program of CSU and Hunan Provincial Innovation Foundation For Postgraduate (CX2018B051).}

\thanks{$^*$ Corresponding author: x.j.wang7@csu.edu.cn; x.j.wang7@gmail.com
             (Xiaojie Wang)}

\begin{document}
\maketitle

\centerline{\scshape Ziheng Chen, Siqing Gan and Xiaojie Wang$^*$}
\medskip
{\footnotesize
 \centerline{School of Mathematics and Statistics, Central South University}
   \centerline{Changsha 410083, Hunan, China}
} 

%
%
%
%

\bigskip


\begin{abstract}
      This paper first establishes a fundamental mean-square convergence theorem for general one-step numerical approximations of L\'{e}vy noise driven  stochastic differential equations with non-globally Lipschitz coefficients.
      Then two novel explicit schemes are designed and their convergence rates are exactly identified via the fundamental theorem.
      Different from existing works, we do not impose a globally Lipschitz condition on the jump coefficient but formulate appropriate assumptions to allow for its super-linear growth. However, we require that the L\'{e}vy measure is finite.
      New arguments are developed to handle essential difficulties in the  convergence analysis, caused by the super-linear growth of the jump coefficient
      and the fact that higher moment bounds of the Poisson increments
      $
          \int_t^{t+h} \int_Z \,\bar{N}(\mbox{d}s,\mbox{d}z),
          t \geq 0, h >0
      $
      contribute to magnitude not more than $O(h)$.
      Numerical results are finally reported to confirm the theoretical findings.
\end{abstract}

\section{Introduction} \label{sec:introduction}

As a class of important mathematical models, stochastic differential equations (SDEs) driven by Gaussian noise have been widely used in
finance, biology, fluid mechanics, chemistry and many other scientific fields. Nevertheless, in the real life one often encounters problems
influenced by event-driven uncertainties, which can be captured by jump component. For instance, the stock price movements might suffer from sudden and significant impacts caused by unpredictable important events such as market crashes, announcements made by central banks, changes in credit ratings, etc. In order to model the event-driven phenomena,  it is necessary and significant to introduce jump-diffusion SDEs, a typical example of non-Gaussian noise (consult, e.g., \cite{Platen10} for more explanation). Since the analytic solutions of nonlinear SDEs with jumps are rarely available, numerical solutions become a powerful tool to understand the behavior of the underlying problems. Therefore this paper concerns the design and mean-square convergence analysis of discrete-time approximations for jump-diffusion SDEs.

Let $d,m \in \N$, $T > 0$ and let $(\Omega,\F,\P)$ be a complete probability space with a filtration $\{\F_{t}\}_{0\leq{t}\leq{T}}$ satisfying the usual conditions.
Let $\{W(t)\}_{0\leq{t}\leq{T}}$ be an $m$-dimensional
$\{\F_{t}\}_{0\leq{t}\leq{T}}$-adapted Wiener process. Let $(Z, \mathcal{Z}, \nu)$ be a measure space with $Z\subseteq\R^d \setminus\{0\}$ and let $N(\mbox{d}t,\mbox{d}z)$ be an $\{\F_{t}\}_{0\leq{t}\leq{T}}$-adapted Poisson random measure defined on $([0,T] \times Z , \mathfrak{B}([0,T] \times Z) )$
with $\nu \neq 0$ and $\nu (Z) < \infty$.
The compensated Poisson random measure is denoted by
$\bar{N}(\mbox{d}t,\mbox{d}z):= N(\mbox{d}t,\mbox{d}z)-\nu(\mbox{d}z)\mbox{d}t$.
We consider the following jump-diffusion SDEs
\begin{equation}\label{eq:JSDE}
      \mbox{d}X(t)
      =
      f(X(t^{-}))\,\mbox{d}t
      +
      g(X(t^{-}))\,\mbox{d}W(t)
      +
      \int_Z\sigma(X(t^{-}),z)\,\bar{N}(\mbox{d}t,\mbox{d}z),
      \enskip
      \forall t\in(0,T]
\end{equation}
with $X(0)=X_{0}$, where $X(t^{-}) = \lim_{s \nearrow t} X(s), \forall t \in (0,T]$ and the drift coefficient $f \colon \R^{d} \to \R^{d}$, the diffusion coefficient $g \colon \R^{d} \to \R^{{d}\times{m}}$ and the jump coefficient $\sigma \colon \R^{d}\times{Z} \to \R^{d}$ are assumed to be deterministic and Borel measurable. Under further assumptions specified later, a unique solution exists in $L^2(\Omega, \R^d)$ for \eqref{eq:JSDE}
and its numerical approximation is a central topic of this work.

As the jump component vanishes, i.e., $\sigma \equiv 0$, the underlying jump-diffusion SDEs reduce to the usual SDEs without jumps, numerical methods for which have been extensively studied for the past decades (consult monographs \cite{Hutzenthaler15,kloeden92numerical,Milstein04} and references therein), in the context of both numerical convergence and numerical stability. Under globally Lipschitz conditions,
the corresponding numerical analysis is well-understood \cite{kloeden92numerical,Milstein04}.
However, coefficients of most models in applications do not obey the classical conditions but, e.g., might behave super-linearly.
Recently, Hutzenthaler, Jentzen and Kloeden \cite{Hutzenthaler10} showed that the standard explicit Euler method produces divergent strong and weak numerical approximations in a finite time interval once one of coefficients grows super-linearly.
By contrast, as already shown by Higham, Mao and Stuart \cite{Higham02}, Mao and Szpruch \cite{Mao13a}, Andersson and Kruse\cite{andersson2016mean}, the backward (implicit) Euler method, computationally much more expensive than the explicit Euler method, can be strongly convergent under certain non-globally Lipschitz conditions.
These observations suggest that special care must be taken to construct and analyze convergent numerical schemes in non-globally Lipschitz setting, and this interesting subject has been investigated in a great portion of the literature
\cite{andersson2016mean,Beyn16,beyn2015stochastic,Dareiotis16,fang2016adaptive,
Hutzenthaler11,Hutzenthaler15,Hutzenthaler10,Hutzenthaler12,hutzenthaler2014perturbation,
hutzenthaler2018exponential,kelly2018adaptiveIMA,kumar2017explicit,
kumar2017tamed,liu2013strong,Mao13b,Mao13a,mao2015truncated,mao2016convergence,
Milstein05,Sabanis16,sabanis2013note,sabanis2019explicit,szpruch2018vintegrability,
Tambue15,Tretyakov13,Wang13,ZM17,Zhang14,zong2014convergence}.
In 2012,  Hutzenthaler, Jentzen and Kloeden \cite{Hutzenthaler12} introduced an explicit method, called the tamed Euler method,
to numerically solve SDEs with super-linearly growing drift coefficients and globally Lipschitz diffusion coefficients.
Since then, various explicit schemes are designed and analyzed for SDEs with (more general) locally Lipschitz coefficients
\cite{Beyn16,beyn2015stochastic,Dareiotis16,Hutzenthaler11,Hutzenthaler15,
hutzenthaler2014perturbation,kumar2017tamed,liu2013strong,mao2015truncated,
mao2016convergence,Sabanis16,sabanis2013note,sabanis2019explicit,
Tambue15,Tretyakov13,Wang13,ZM17,Zhang14,zong2014convergence}.
Readers can, e.g.,  refer to \cite{Hutzenthaler15}
for a more comprehensive list of references.
Particularly, we should mention a closely relevant article \cite{Tretyakov13} by Tretyakov and Zhang, where a fundamental strong convergence theorem was derived
in a non-globally Lipschitz setting, giving an extension of a counterpart in the globally Lipschitz setting \cite{Milstein87,Milstein04}.
Moreover, an explicit balanced Euler method, given by
\begin{equation*}\label{eq:intro-Tretyakov-Zhang-Euler}
      Y_{n+1}
      =
      Y_{n}
      +
      \frac{f(Y_n)h + g(Y_n)\Delta{W_n}}
           {1+|f(Y_n)|h+|g(Y_n)\Delta{W_n}|},
      \quad
      Y_{0} = X_{0}
\end{equation*}
and a fully implicit Euler method are examined and their convergence rates are obtained, with the aid of the fundamental convergence theorem.
Another two closely related papers are \cite{ZM17,Zhang14} by Zhang and Ma, introducing a sine Euler method, defined by
\begin{equation*}\label{eq:intro-sine-Euler-Zhang}
      Y_{n+1}
      =
      Y_{n}
      +
      \sin( f(Y_n)h )
      +
      \sin( g(Y_n)\Delta{W_n} ),
      \quad
      Y_{0} = X_{0}.
\end{equation*}
%

When $\sigma \not \equiv 0$, the underlying jump-diffusion problem, as a typical non-continuous stochastic process, has been increasingly studied in recent years and a lot of progress has been achieved on numerical analysis of explicit and implicit schemes \cite{Bruti07,Deng19, Gardon04,Higham05,Higham07,hu2011convergence,
Jacod05,Kohatsu10,Liu00,Maghsoodi96,Platen10,Wang10,yang2017transformed}.
Particularly, some explicit time-stepping schemes are very recently introduced and their convergence rates are analyzed in non-globally Lipschitz setting \cite{Dareiotis16,Deng19,kumar2017explicit,kumar2017tamed,Tambue15}.
However, all existing works on convergence of numerical methods for SDEs with jumps, to the best of our knowledge,
impose globally Lipschitz conditions on the jump coefficient (consult the very recent publications \cite{Dareiotis16,kumar2017explicit,kumar2017tamed}
and references therein).
As pointed out in Chapter 1, Section 9, on Page 59 of \cite{Platen10} by Platen and Bruti-Liberati,
for certain applications, the Lipschitz condition on the jump coefficient is too restrictive.
For instance, for modeling state-dependent intensities, as discussed in Sect. 1.8 therein,
it is convenient to use jump coefficients that are not Lipschitz continuous.
This indicates that SDEs with non-globally Lipschitz continuous jump coefficients have applications in certain fields
and motivates the present numerical analysis in a more general setting, allowing for non-globally Lipschitz continuous jump coefficients.
%
%

We first establish a fundamental mean-square convergence theorem for general one-step numerical methods under certain non-globally Lipschitz conditions (see Assumptions \ref{as:assumption1}, \ref{as:assumption2}).
%
Although the proof of the fundamental theorem follows the basic lines in previous works \cite{Milstein87,Milstein04,Tretyakov13}, some extension of their arguments are made due to the presence of the jump term.  For example, new techniques are required and new assumptions (see Assumptions \ref{as:assumption2} and \ref{ass:drift-polynomial}) are formulated,  to treat additional terms resulting from a jump version of the It\^{o} formula.
As applications of the fundamental theorem, a new version of the tamed Euler method
\begin{equation}\label{eq:intro.tame.method}
      Y_{n+1}
      =
      Y_{n}
      +
      \frac{f(Y_n)h}{1+|f(Y_n)|h}
      +
      \frac{g(Y_n)\Delta{W_n}}{1+|g(Y_n)|h}
      +
      \int_{t_n}^{t_{n+1}}\hspace{-0.75em}\int_{Z}
          \frac{\sigma(Y_n,z)}{1+|\sigma(Y_n,z)|h}
      \,\bar{N}(\mbox{d}s,\mbox{d}z)
\end{equation}
and a so-called sine Euler method
\begin{equation}\label{eq:intro.sin}
      Y_{n+1}
      =
      Y_n
      +
      \sin(f(Y_n)h)
      +
      \frac{\sin(g(Y_n)h)}{h}\Delta{W_n}
      +
      \int_{t_n}^{t_{n+1}}\hspace{-0.75em}\int_Z
          \frac{\sin(\sigma(Y_n,z)h)}{h}
      \,\bar{N}(\mbox{d}s,\mbox{d}z)
\end{equation}
are carefully constructed and their mean-square convergence rates
are accordingly identified (see Theorems \ref{Theorem:tame.order.I}
and \ref{thm:sinerror.sub.nonadditive}).
The most challenging and technical part in the applications of the fundamental theorem lies on proving boundedness of higher order moments of numerical approximations (see Subsections \ref{subset:bounded-moments} and \ref{subsec:bonded-mom}). Unlike the Wiener increments $W(t+h)-W(t), t \geq 0, h >0$, higher moment bounds of the Poisson increments
$\int_t^{t+h} \int_Z \,\bar{N}(\mbox{d}s,\mbox{d}z), t \geq 0,h >0$
contribute to magnitude not more than $O(h)$. This significant difference, together with the possible super-linear growth of the jump coefficient,
makes the approach used in \cite{Tretyakov13} unworkable here since the nice property of Wiener increments
$\E [ \| W(t+h)-W(t) \|^l ] = O(h^{\frac{l}{2}}), l \in \N $ was essentially used there (see the treatment of the last term of (3.6) in \cite{Tretyakov13}), where $\|\cdot\|$ denotes the Euclidean vector norm in $\R^{m}$.
This forces us to develop new arguments for the present jump setting.
In short, we work with continuous-time approximations of \eqref{eq:intro.tame.method} and \eqref{eq:intro.sin} and carry out rather careful and delicate estimates for all involved terms (see Remark \ref{remark:Mom-Bound} and the proof of Lemma \ref{lemma:tame.moment.bound.a}).
%
Equipped with bounded numerical moments, we examine the local truncation errors of the schemes.
This together with the fundamental theorem helps us to obtain the mean-square convergence rates arbitrarily close to the classical order $\tfrac12$
(see Theorems \ref{Theorem:tame.order.I} and \ref{thm:sinerror.sub.nonadditive}).
To the best of our knowledge, this is the first result to identify convergence rates of numerical approximations of jump-diffusion SDEs with possibly super-linearly growing jump coefficients. Moreover, when the jump component vanishes, i.e., $\sigma \equiv 0$,
an exact order $\tfrac12$ can be attained, see Corollaries \ref{Theorem:nonjump} and \ref{thm:sinerror.sub.nonadditive.sigmazero}, which recovers the relevant results in \cite{andersson2016mean,Hutzenthaler15,Hutzenthaler12,Mao13a,
Sabanis16,sabanis2013note,Tretyakov13}.

Now  we compare convergence results in this article with corresponding results in existing literature.
The contributions \cite{Dareiotis16,kumar2017explicit}  derived the strong convergence rates of different tamed Euler methods
for a wider class of L\'{e}vy SDEs by allowing $\nu(Z) = \infty$, but
with the jump coefficients satisfying the globally Lipschitz conditions, see Assumptions B-3 in \cite{Dareiotis16} and A-8 in \cite{kumar2017explicit}. In contrast,
we reformulate a more relaxed condition on the jump coefficient to allow for its super-linear growth.
Here we require $\nu(Z)<\infty$, which is a limitation compared with \cite{Dareiotis16,kumar2017explicit} and only covers
SDEs driven by Wiener process and compound Poisson process.
Furthermore, in \cite{Dareiotis16,kumar2017explicit} the $\mathcal{L}^{p}$-convergence rates (for any $p \geq 2$) were obtained,
but here we can only get mean-square convergence rates.
Finally, to achieve bounded moments of the exact solution, we first show bounded $\bar{p}$-th moments for
the highest order $\bar{p}$ being a sufficiently large even number and the general $p$-th moments with arbitrary $p \leq \bar{p}$
follow immediately from the H\"{o}lder inequality, see Theorem \ref{th:exactbound}.

The remainder of this paper is organized as follows. The next section concerns properties of jump-diffusion SDEs. A fundamental mean-square convergence theorem for general one-step approximations is established
in Section \ref{sec:theorem}. In Sections \ref{sec:tamed} and \ref{sec:sine}, we propose two new explicit schemes and identify their mean-square convergence rates via the obtained fundamental theorem. Finally numerical experiments are performed to illustrate the theoretical results.

\section{L\'{e}vy noise driven SDEs with non-globally Lipschitz coefficients}
\label{sec:problem}

We start with some notation. Let $|\cdot|$ and $\left<\cdot,\cdot\right>$ be the Euclidean norm and the inner product in $\R^d$. By $A^{T}$ we denote the transpose of a vector or matrix $A$. If $A$ is a matrix, we let $|A|=\sqrt{\tr(A^{T}A)}$ be its trace norm. If $B$ is a set, its complement and indicator function are denoted by $B^c$ and $\1_B$, respectively. We denote by $L^{r}_{\F_{0}} (\Omega, \R^d ), r \in \N $ the family of $\R^d$-valued $\F_{0}$-measurable random variables $\xi$ with $\E[|\xi|^{r}]<\infty$. For simplicity, the letter $K$ denotes a generic positive constant that is independent of time stepsize and varies for each appearance.

An $\{\F_{t}\}_{0 \leq {t} \leq {T}}$-adapted $\R^d$-valued stochastic process $\{X(t)\}_{0\leq{t}\leq{T}}$ is called a solution of \eqref{eq:JSDE} if it is almost surely right continuous with left limits and satisfies
\begin{equation}\label{eq:SIDE}
\begin{split}
      X(t)
      =
      X_0&+\int_0^tf(X(s^-))\,\mbox{d}s
      +
      \int_0^tg(X(s^-))\,\mbox{d}W(s)
      \\&+
      \int_0^t\int_Z\sigma(X(s^-),z)\,\bar{N}(\mbox{d}s,\mbox{d}z),
      ~~ \P\text{-a.s.}, ~~ \forall t \in [0,T].
\end{split}
\end{equation}
Let us make the following assumption.

\begin{assumption}\label{as:assumption1}
      Let $X_0 \in L^{2}_{\F_{0}} (\Omega, \R^d )$ and let the mappings $f \colon \R^{d} \to \R^{d}$, $g \colon \R^{d} \to \R^{{d}\times{m}}$, $\sigma(\cdot,z) \colon \R^{d} \to \R^{d},~\forall z\in{}Z$ be continuous and satisfy the monotone condition, i.e., there exists $K>0$ such that for all $x,y\in\R^d$,
      \begin{equation}\label{eq:gmc}
            2\langle x-y,f(x)-f(y) \rangle
            +
            |g(x)-g(y)|^{2}
            +
            \int_{Z}|\sigma(x,z)-\sigma(y,z)|^{2}\,\nu(\mbox{d}z)
            \leq
            K|x-y|^2.
      \end{equation}
Also let the coefficients satisfy the coercivity condition, i.e., there exists $K > 0$ such that
\begin{equation}\label{eq:gcc.a}
      2\langle x,f(x) \rangle
      +
      |g(x)|^{2}
      +\int_{Z}|\sigma(x,z)|^2\,\nu(\mbox{d}z)
      \leq
      K(1+|x|^2),
      \quad
      \forall x \in \R^d.
\end{equation}
\end{assumption}

Let us give the existence of the unique solution to \eqref{eq:JSDE},
see \cite[Theorem 1]{Gyongy80}.

\begin{theorem}
       Let Assumption \ref{as:assumption1} be satisfied. Then \eqref{eq:JSDE} admits a unique solution $\{X(t)\}_{0\leq{t}\leq{T}}$ given by \eqref{eq:SIDE} and almost all its trajectories are right-continuous with left limits on $[0,T]$.
\end{theorem}

To discuss the higher order moments of $\{X(t)\}_{0\leq{t}\leq{T}}$, we make the following assumption.

\begin{assumption}\label{as:assumption2}

      Let $\varepsilon$ be an arbitrary positive number and let $\bar{p}\geq2$ be a sufficiently large even number. Assume that $X_{0} \in L^{\bar{p}}_{\F_{0}} (\Omega, \R^d )$ and that there exists $K > 0$ such that
      \begin{equation}\label{eq:gcc.b}
      \begin{split}
            \bar{p}|x|^{\bar{p}-2}&\langle x,f(x) \rangle
            +
            \tfrac{\bar{p}(\bar{p}-1)}{2}|x|^{\bar{p}-2}|g(x)|^{2}
            \\&+
            \big(1+(\bar{p}-2)\varepsilon\big)
            \int_Z|\sigma(x,z)|^{\bar{p}}\,\nu(\mbox{d}z)
            \leq
            K(1+|x|^{\bar{p}}),
            \quad  \forall x \in \R^{d}.
      \end{split}
      \end{equation}
\end{assumption}

Note that the assumption \eqref{eq:gcc.b} with $\bar{p} = 2$ reduces to \eqref{eq:gcc.a}. We also point out that $\bar{p}$ is the upper order of the bounded moments of the exact solution (see Theorem \ref{th:exactbound})
and the parameter $\varepsilon$ comes from proving bounded moments of the exact solution up to order $\bar{p}$
(see \eqref{eq:analbound.inequality} and \eqref{eq:analbound.gcc} below). With the above setting, we arrive at the following result.

\begin{theorem}\label{th:exactbound}
      Suppose Assumptions \ref{as:assumption1}, \ref{as:assumption2} hold and let $\{X(t)\}_{0\leq{t}\leq{T}}$ be defined by \eqref{eq:SIDE}. Let $\bar{p}\geq2$ coming from \eqref{eq:gcc.b} be a sufficiently large even number. Then there exists $K>0$ such that
      \begin{equation}\label{eq:anal.solu.bounded}
            \E\big[|X(t)|^{\bar{p}}\big]
            \leq
            K\big(1+\E\big[|X_{0}|^{\bar{p}}\big]\big),
            \quad \forall  t \in [0,T].
      \end{equation}
\end{theorem}

Before presenting its proof, we provide two elementary facts as follows.

\begin{lemma}\label{lemma:inequality}
      Let $\rho \in \N$ and $\varepsilon >0$, it holds that
      \begin{equation}\label{eq:inequalitybino}
            |x+y|^{2\rho} - |x|^{2\rho} - 2\rho\langle x,y \rangle|x|^{2\rho-2}
            \leq
            K|x|^{2\rho} + \big(1+(2\rho-2)\varepsilon\big)|y|^{2\rho},
            \quad
            \forall x,y \in \R^{d}.
      \end{equation}
\end{lemma}
\begin{proof}
It follows from the binomial formula that
\begin{equation}\label{eq:KKK}
\begin{split}
      |x+y|^{2\rho}
      =&
      \big(|x|^{2}+2\left<x,y\right>+|y|^{2}\big)^{\rho}
      =
      \sum\limits_{i=0}^{\rho}\mbox{C}_{\rho}^{i}
      |y|^{2i}\big(2\left<x,y\right>+|x|^{2}\big)^{\rho-i}
      \\=&
      |y|^{2\rho}
      +
      \sum\limits_{j=0}^{\rho}\mbox{C}_{\rho}^{j}
      \big(2\left<x,y\right>\big)^{j}|x|^{2(\rho-j)}
      +
      \sum\limits_{i=1}^{\rho-1}\mbox{C}_{\rho}^{i}
      |y|^{2i}\big(2\left<x,y\right>+|x|^{2}\big)^{\rho-i}
      \\=&
      |y|^{2\rho}+|x|^{2\rho}+2\rho\left<x,y\right>|x|^{2\rho-2}
      +
      \sum\limits_{j=2}^{\rho}\mbox{C}_{\rho}^{j}
      \big(2\left<x,y\right>\big)^{j}|x|^{2(\rho-j)}
      \\&+
      \sum\limits_{i=1}^{\rho-1}\mbox{C}_{\rho}^{i}
      |y|^{2i}\big(2\left<x,y\right>+|x|^{2}\big)^{\rho-i}.
\end{split}
\end{equation}
Employing the Schwarz inequality and the weighted Young inequality guarantees
\begin{equation}\label{eq:K1}
      \mbox{C}_{\rho}^{j}
      \big(2\left<x,y\right>\big)^{j}|x|^{2(\rho-j)}
      \leq
      2^{j}\mbox{C}_{\rho}^{j} |x|^{2\rho-j}|y|^{j}
      \leq
      K_{1}(j)|x|^{2\rho} + \varepsilon|y|^{2\rho},
      \enskip
      j=2,\ldots,\rho,
\end{equation}
where
\begin{equation*}
      K_{1}(j)
      =
      \tfrac{2\rho-j}{2\rho}
      \left( 2^{j}\mbox{C}_{\rho}^{j} \right)^{\frac{2\rho}{2\rho-j}}
      \left( \tfrac{2\rho\varepsilon}{j} \right)^{\frac{j}{j-2\rho}},
      \enskip
      j=2,\ldots,\rho.
\end{equation*}
Applying the techniques used in \eqref{eq:K1} and the inequality
\begin{equation}\label{eq:elementaryinequality}
  \big(\sum\limits_{i = 1}^{n} |a_{i}|\big)^{p}
  \leq
  n^{p} \sum\limits_{i = 1}^{n} |a_{i}|^{p},
  \quad \forall p > 0, a_{i} \in \R, i = 1,2, \ldots, n, n \in \N,
\end{equation}
one gets
\begin{equation}\label{eq:K2K3}
\begin{split}
      \mbox{C}_{\rho}^{i}|y|^{2i}
      \big( 2\left<x,y\right> + |x|^{2} \big)^{\rho-i}
      \leq&
      2^{2\rho-2i}\mbox{C}_{\rho}^{i}|x|^{\rho-i}|y|^{\rho+i}
      +
      2^{\rho-i}\mbox{C}_{\rho}^{i}|x|^{2\rho-2i}|y|^{2i}
      \\\leq&
      \big( K_2(i)+K_3(i) \big)|x|^{2\rho}
      +
      \varepsilon|y|^{2\rho},
\end{split}
\end{equation}
where
$$
    K_2(i)
    =
    \frac{\rho-i}{2\rho}
    \big( 2^{2\rho-2i}\mbox{C}_{\rho}^{i} \big)^{\frac{2\rho}{\rho-i}}
    \big( \frac{\rho\varepsilon}{\rho+i} \big)^{\frac{\rho+i}{i-\rho}},
    ~~
    K_3(i)
    =
    \frac{\rho-i}{\rho}
    \big( 2^{\rho-i}\mbox{C}_{\rho}^{i} \big)^{\frac{\rho}{\rho-i}}
    \big( \frac{\rho\varepsilon}{2i} \big)^{\frac{i}{i-\rho}}
$$
for all $i=1,\ldots,\rho-1$. Then \eqref{eq:inequalitybino} follows by inserting \eqref{eq:K1} and \eqref{eq:K2K3} into \eqref{eq:KKK}.
\end{proof}

The next lemma is the It\^{o} formula
\cite[Theorem 33]{Protter90}, frequently used later.
\begin{lemma}\label{lemma:itoformula}
      Let $\{X(t)\}_{t \geq 0}$ be a $\R^{d}$-valued stochastic process characterized by \eqref{eq:JSDE} and let $V \colon \R^{d} \to \R$ be a continuously twice differentiable function. Then $V(X(t))$ is given by
      \begin{equation*}
      \begin{split}
            \mbox{d}V(X(t))
            =&
            [
              V_{x}(X(t^{-}))f(X(t^{-}))
              +
              \tfrac{1}{2}
              \tr( g^{T}(X(t^{-})) V_{xx}(X(t^{-}) g(X(t^{-})) )
            ] \mbox{d}t
            \\&+
            \int_{Z} V(X(t^{-})+\sigma(X(t^{-}),z))-V(X(t^{-})) \bar{N}(\mbox{d}t,\mbox{d}z)
            \\&+
            V_{x}(X(t^{-}))g(X(t^{-}))\,\mbox{d}W(t)
            +
            \int_{Z} V(X(t^{-})+\sigma(X(t^{-}),z))
            \\&-V(X(t^{-}))-V_{x}(X(t^{-}))\sigma(X(t^{-}),z)            \nu(\mbox{d}z)\mbox{d}t,
            \quad \P\text{-a.s.}, \quad \forall t \geq 0.
      \end{split}
      \end{equation*}
\end{lemma}


At the moment, we are ready to prove Theorem \ref{th:exactbound}.

\begin{proof}[Proof of Theorem \ref{th:exactbound}.]
For every integer $k\geq1$, define the stopping times
$\tau_k = \inf\{t\in[0,T]:|X(t)|>k\}, k \in \N$.
Clearly $\tau_k \uparrow {T}$ as $k \rightarrow \infty, \P\text{-a.s.}$ and
$|X(t)| \leq {k}$ 
for all $0 \leq {t} \leq \tau_k$.
By Lemma \ref{lemma:itoformula} we can derive that for all $t \in [0,T]$, 
\begin{equation*}\label{eq:analbound.Ito}
\begin{split}
      |X(t)|^{\bar{p}}
      =&
      |X_0|^{\bar{p}}
      +
      {\bar{p}}\int_0^t
           |X(s^-)|^{{\bar{p}}-2}
           \langle X(s^-),f(X(s^-)) \rangle
       \,\mbox{d}s
      \\&+
      \frac{{\bar{p}}}{2}
      \int_0^t
          |X(s^-)|^{{\bar{p}}-2}|g(X(s^-))|^2
      \,\mbox{d}s
      \\&+
      \frac{{\bar{p}}({\bar{p}}-2)}{2}
      \int_0^t
          |X(s^-)|^{{\bar{p}}-4}|X(s^-)^{T}g(X(s^-))|^2
      \,\mbox{d}s
      \\&+
      {\bar{p}}\int_0^t
          |X(s^-)|^{{\bar{p}}-2}
          \langle X(s^-),g(X(s^-)) \,\mbox{d}W(s) \rangle
      \\&+
      \int_0^t\int_Z
          |X(s^-)+\sigma(X(s^-),z)|^{\bar{p}} - |X(s^-)|^{\bar{p}}
      \,N(\mbox{d}s,\mbox{d}z)
      \\&-
      {\bar{p}}\int_0^t\int_Z
           |X(s^-)|^{{\bar{p}}-2}
           \langle X(s^-),\sigma(X(s^-),z) \rangle
       \,\nu(\mbox{d}z)\mbox{d}s, ~~ \P\text{-a.s.}.
\end{split}
\end{equation*}
This together with Lemma \ref{lemma:inequality} results in
\begin{equation}\label{eq:analbound.inequality}
\begin{split}
      |X(t)|^{\bar{p}}
      \leq&
      |X_0|^{\bar{p}}
      +
      {\bar{p}} \! \int_0^t
           |X(s^-)|^{{\bar{p}}-2}
           \langle X(s^-),f(X(s^-))\rangle
       \,\mbox{d}s
      \\&
      +
      \frac{{\bar{p}}({\bar{p}}-1)}{2}
      \!
      \int_0^t
          |X(s^-)|^{{\bar{p}}-2}|g(X(s^-))|^2
      \,\mbox{d}s
      \\&+
      {\bar{p}}\int_0^t
           |X(s^-)|^{{\bar{p}}-2}
           \langle X(s^-),g(X(s^-)) \,\mbox{d}W(s) \rangle
      \\&+
      \int_0^t\int_Z
          |X(s^-)+\sigma(X(s^-),z)|^{\bar{p}}
          -
          |X(s^-)|^{\bar{p}}
      \,\bar{N}(\mbox{d}s,\mbox{d}z)
      \\&+
      \big(1+({\bar{p}}-2)\varepsilon\big)
      \int_0^t\int_Z
          |\sigma(X(s^-),z)|^{\bar{p}}
      \,\nu(\mbox{d}z)\mbox{d}s
      \\&+
      K\nu(Z)\int_0^t |X(s^-)|^{\bar{p}}\,\mbox{d}s,
       ~~ \P\text{-a.s.}.
\end{split}
\end{equation}
Thanks to \eqref{eq:gcc.b}, $\nu(Z) < \infty$ and the martingale property, we deduce
\begin{equation}\label{eq:analbound.gcc}
\begin{split}
      \E\big[  |X(t\wedge\tau_k)|^{\bar{p}}  \big]
      \leq&
      K
      +
      \E\big[|X_0|^{\bar{p}}\big]
      +
      K\int_{0}^{t}
           \E\big[ |X(s^-\wedge\tau_k)|^{\bar{p}} \big]
       \,\mbox{d}s
      \\\leq&
      K
      +
      \E\big[|X_0|^{\bar{p}}\big]
      +
      K\int_{0}^{t}
           \sup\limits_{0 \leq r \leq s}\E\big[ |X(r\wedge\tau_k)|^{\bar{p}} \big]
       \,\mbox{d}s,
\end{split}
\end{equation}
where the properties of the right-continuous with left limits functions \cite[Section 2.9]{applebaum2009levy} were used in the last step.
This immediately gives
\begin{equation}\label{eq:analbound.gcccc}
      \sup\limits_{0 \leq r \leq t}
      \E\big[  |X(r\wedge\tau_k)|^{\bar{p}}  \big]
      \leq
      K
      +
      \E\big[ |X_0|^{\bar{p}} \big]
      +
      K\int_{0}^{t}
           \sup\limits_{0 \leq r \leq s}\E\big[ |X(r\wedge\tau_k)|^{\bar{p}} \big]
       \,\mbox{d}s.
\end{equation}
By the Gronwall inequality, one gets
\begin{equation*}\label{eq:analbound.gronwall}
      \sup\limits_{0 \leq r \leq t}
      \E\big[ |X(r \wedge \tau_k)|^{\bar{p}} \big]
      \leq
      \big( K + \E\big[ |X_{0}|^{\bar{p}} \big] \big)e^{KT},
      \quad \forall t \in [0,T].
\end{equation*}
Letting $t = T$, $k\rightarrow\infty$ and applying Fatou's lemma finish the proof.
\end{proof}

Before ending this section, we would like to point out that the condition $\nu (Z) < \infty$ can not be relaxed following our approach. However, one can consult \cite{Dareiotis16, kumar2017explicit}, treating SDEs with globally Lipschitz jump coefficients, where the authors used a different approach for the convergence analysis and $ \nu (Z) = \infty $ was allowed.

\section{A fundamental mean-square convergence theorem}
\label{sec:theorem}

This section aims to establish a fundamental mean-square convergence theorem
for general one-step approximation. Given a uniform mesh on $[0,T]$ with $h = \tfrac{T}{N}, N \in \N$ being the stepsize, we denote by $X_{t,x}(t+h)$ the solution of \eqref{eq:JSDE} at $t+h$ with initial value $X_{t,x}(t) = x$. For $x \in \R^{d}$, $t \in [0,T]$, $h >0$, $0 < t+h \leq T$, the general one-step approximation $Y_{t,x}(t+h)$ for $X_{t,x}(t+h)$, depends on $t, x, h, W(t+h)-W(t), \int_t^{t+h} \int_Z \,\bar{N}(\mbox{d}s,\mbox{d}z)$ and is given by
\begin{equation*}\label{eq:one.step}
      Y_{t,x}(t+h)
      =
      x
      +
      \Psi\Big( t, x, h,
                W(t+h)-W(t),
                \int_t^{t+h} \int_Z \,\bar{N}(\mbox{d}s,\mbox{d}z)
          \Big),
\end{equation*}
where $\Psi$ is a function from $[0,T] \times R^{d} \times (0,T] \times \R^{m} \times \R^{d}$ to $\R^{d}$. By $Y_{t,x}( t +h)$ we denote an approximation of the solution at $t+h$ with initial value $Y_{t,x}( t) = x $. Then one can use $Y_{n+1} = Y_{t_n,Y_n}(t_{n+1})$
to recurrently construct numerical approximations $\{Y_n\}_{0\leq{n}\leq{N}}$ on the uniform mesh grid $\{ t_n = n h, n= 0,1,\ldots,N \}$, given by
\begin{equation}\label{eq:one.step.numer}
      Y_{n+1}
      =
      Y_{n}
      +
      \Psi\Big(t_n,Y_n,h,
      W(t_n+h)-W(t_n), \int_{t_{n}}^{t_{n}+h}\hspace{-0.75em}\int_Z \,\bar{N}(\mbox{d}s,\mbox{d}z)
      \Big)
\end{equation}
for all $n=0, 1, \ldots,N-1$ with $Y_0=X_0$. Alternatively, one can write
\begin{equation*}
      Y_{n+1}
      =
      Y_{t_n,Y_n}(t_{n+1})
      =
      Y_{t_0,Y_0}(t_{n+1}),
      \quad
      n=0, 1, \ldots, N-1.
\end{equation*}

To proceed, we need the following assumption.

\begin{assumption}\label{ass:drift-polynomial}
      Assume the drift coefficient $f$ of \eqref{eq:JSDE} behaves polynomially in the sense that there exist $K,q \geq 0$ such that
      \begin{equation}\label{eq:fpg}
            |f(x)-f(y)|^{2}
            \leq
            K(1+|x|^{q}+|y|^{q})|x-y|^{2},
            \quad \forall x,y \in \R^{d}.
     \end{equation}
\end{assumption}
The inequality \eqref{eq:fpg} implies the polynomial growth bound
\begin{equation}\label{eq:fsg}
      |f(x)|^2 \leq K(1+|x|^{2+q}),
      \quad \forall x \in \R^{d}.
\end{equation}
This together with \eqref{eq:gcc.b} further shows that 
\begin{align}
      \label{eq:gsg}
      |x|^{\bar{p}-2}|g(x)|^2
      \leq&
      K(1+|x|^{\bar{p}+\frac{q}{2}}),
      \quad \forall x \in \R^{d},
      \\
      \label{eq:hsg}
      \int_Z |\sigma(x,z)|^{\bar{p}} \,\nu(\mbox{d}z)
      \leq&
      K(1+|x|^{\bar{p}+\frac{q}{2}}),
      \quad \forall x \in \R^{d}.
\end{align}
Moreover, by virtue of \eqref{eq:gmc} and \eqref{eq:fpg}, it is easy to derive
\begin{align}
      \label{eq:gpg}
      |g(x)-g(y)|^{2}
      \leq&
      K(1+|x|^{\frac{q}{2}}+|y|^{\frac{q}{2}})|x-y|^{2},
      \quad \forall x,y \in \R^{d},
      \\\label{eq:hpg}
      \int_Z|\sigma(x,z)-\sigma(y,z)|^{2}\,\nu(\mbox{d}z)
      \leq&
      K(1+|x|^{\frac{q}{2}}+|y|^{\frac{q}{2}})|x-y|^{2},
      \quad \forall x,y \in \R^{d}.
\end{align}

Here and below we always assume $h \in (0,h_{0}]$ for some $h_{0} \in (0,1)$.

\begin{lemma}\label{LemmaA}
      Suppose Assumptions \ref{as:assumption1}, \ref{as:assumption2}, \ref{ass:drift-polynomial} hold. Let $\bar{p}$ and $q$ come from \eqref{eq:gcc.b} and \eqref{eq:fpg}, respectively, satisfying $q \leq \bar{p}$. For the representation
      \begin{equation}\label{eq:repre}
            X_{t,x}(t+h)-X_{t,y}(t+h)
            =
            x-y+\Phi_{t,x,y}(t+h),
            \quad \forall x,y \in \R^d,
      \end{equation}
      there exists
      $K > 0$ depending on $h_0 \in (0, 1)$ but independent of $h$  such that
      for all $0 < h \leq h_{0}, t+h \leq T$,
      \begin{align}
            \label{eq:tx}
            \E\big[ |X_{t,x}(t+h)-X_{t,y}(t+h)|^2 \big]
            \leq
            (1+Kh)|x-y|^2,
            \quad &\forall x,y \in \R^d,
            \\
            \label{eq:txy}
            \E\big[ |\Phi_{t,x,y}(t+h)|^2 \big]
            \leq
            Kh(1+|x|^q+|y|^q)^{\frac{1}{2}}|x-y|^2,
            \quad &\forall x,y \in \R^d.
      \end{align}
\end{lemma}
\begin{proof}
By Lemma \ref{lemma:itoformula}, we infer that for all $0\leq\theta\leq{h}$,
\begin{equation*}
\begin{split}
      &\big|X_{t,x}(t+\theta)-X_{t,y}(t+\theta)\big|^{2}
      \\=&
      |x-y|^2
      +
      \int_{t}^{t+\theta}
          |g(X_{t,x}(s^{-}))-g(X_{t,y}(s^-))|^2
      \,\mbox{d}s
      \\
      \end{split}
\end{equation*}\begin{equation*}\label{eq:tx.uIto}
\begin{split}&+
      2\int_{t}^{t+\theta}
           \langle X_{t,x}(s^-)-X_{t,y}(s^-),f(X_{t,x}(s^-))-f(X_{t,y}(s^-))
           \rangle
       \,\mbox{d}s
      \\&-
      2\int_{t}^{t+\theta}\hspace{-0.5em}\int_{Z}
           \langle
           X_{t,x}(s^-)-X_{t,y}(s^-)
           ,
           \sigma(X_{t,x}(s^-),z)-\sigma(X_{t,y}(s^-),z)
           \rangle
       \,\nu(\mbox{d}z)\mbox{d}s
      \\&+
      2\int_t^{t+\theta}
           \langle
           X_{t,x}(s^-)-X_{t,y}(s^-)
           ,
           (g(X_{t,x}(s^-))-g(X_{t,y}(s^-)))\,\mbox{d}W(s)
           \rangle
      \\&+
      \int_t^{t+\theta}\hspace{-0.5em}\int_{Z}
          |
            X_{t,x}(s^-)-X_{t,y}(s^-)
            +
            \sigma(X_{t,x}(s^-),z)-\sigma(X_{t,y}(s^-),z)
          |^2
      \\&-
          |  X_{t,x}(s^-) - X_{t,y}(s^-)  |^2
      \,N(\mbox{d}s,\mbox{d}z), ~~ \P\text{-a.s.}.
\end{split}
\end{equation*}
Taking expectations
and applying the martingale property imply
\begin{equation*}\label{eq:tx.umartingale}
\begin{split}
      &\E\big[
        |X_{t,x}(t+\theta)-X_{t,y}(t+\theta)|^2
        \big]
      \\=&
      |x-y|^2
      +
      \int_{t}^{t+\theta}
          \E\big[ |g(X_{t,x}(s^{-}))-g(X_{t,y}(s^-)) |^2\big]
      \,\mbox{d}s
      \\&+
      2\int_{t}^{t+\theta}
           \E\big[
             \langle
             X_{t,x}(s^-)-X_{t,y}(s^-)
             ,
             f(X_{t,x}(s^-))-f(X_{t,y}(s^-))
             \rangle
             \big]
       \,\mbox{d}s
      \\&+
      \int_t^{t+\theta}
          \E\Big[
            \int_Z
            |\sigma(X_{t,x}(s^-),z)-\sigma(X_{t,y}(s^-),z)|^2
            \,\nu(\mbox{d}z)
            \Big]
      \mbox{d}s.
\end{split}
\end{equation*}
By the techniques used in \eqref{eq:analbound.gcc}--\eqref{eq:analbound.gcccc}, we derive from \eqref{eq:gmc} that
\begin{equation*}
\begin{split}
      &\sup\limits_{t \leq r \leq t+\theta}
      \E\big[ |X_{t,x}(r)-X_{t,y}(r)|^2 \big]
      \\\leq&
      |x-y|^2
      +
      K\int_t^{t+\theta}
           \sup\limits_{t \leq r \leq s}
           \E\big[|X_{t,x}(r)-X_{t,y}(r)|^2\big]
       \,\mbox{d}s.
\end{split}
\end{equation*}
The Gronwall inequality and the assumption $h \in (0,h_{0}]$ yield
\begin{equation}\label{eq:tx.ugronw}
      \sup\limits_{t \leq r \leq t+\theta}
      \E\big[ |X_{t,x}(r)-X_{t,y}(r)|^2 \big]
      \leq
      \mbox{e}^{Kh}|x-y|^2
      \leq
      (1+\bar{K}h)|x-y|^2
\end{equation}
due to $e^{Kh} \leq 1+\bar{K}h$, where $\bar{K}$ is independent of $h$ but dependent on $h_{0}$. Then \eqref{eq:tx} is validated by letting $\theta=h$ in \eqref{eq:tx.ugronw}.
For \eqref{eq:txy}, using Lemma \ref{lemma:itoformula} shows
\begin{equation*}\label{eq:3.17}
\begin{split}
      |\Phi_{t,x,y}(t&+\theta)|^2
      =
      \int_{t}^{t+\theta}
          |g(X_{t,x}(s^{-}))-g(X_{t,y}(s^-))|^2
      \,\mbox{d}s
      \\&+
      2\int_{t}^{t+\theta}
           \langle
           \Phi_{t,x,y}(s^-),f(X_{t,x}(s^-))-f(X_{t,y}(s^-))
           \rangle
       \,\mbox{d}s
      \\&-
      2\int_{t}^{t+\theta}\int_Z
           \langle
           \Phi_{t,x,y}(s^-),\sigma(X_{t,x}(s^-),z)-\sigma(X_{t,y}(s^-),z)
           \rangle
       \,\nu(\mbox{d}z)\mbox{d}s
      \\&+
      2\int_t^{t+\theta}
           \langle
           \Phi_{t,x,y}(s^-)
           ,
           \left(g(X_{t,x}(s^-))-g(X_{t,y}(s^-))\right)\mbox{d}W(s)
           \rangle
      \\&+
      \int_t^{t+\theta}\int_Z
          |
            \Phi_{t,x,y}(s^-)
            +
            \sigma(X_{t,x}(s^-),z)-\sigma(X_{t,y}(s^-),z)
          |^2
      \\&-
          | \Phi_{t,x,y}(s^-) |^2
      \,N(\mbox{d}s,\mbox{d}z),
      ~ \P\text{-a.s.}.
\end{split}
\end{equation*}
Taking expectations and using \eqref{eq:repre} lead to
\begin{equation*}
\begin{split}
      \E\big[|\Phi_{t,x,y}(t&+\theta)|^2\big]
      =
      \int_{t}^{t+\theta}
           \E\big[ |g(X_{t,x}(s^{-}))-g(X_{t,y}(s^-))|^2 \big]
      \,\mbox{d}s
      \\&+
      2\int_{t}^{t+\theta}
           \E\big[
             \langle
             X_{t,x}(s^-)-X_{t,y}(s^-)
             ,
             f(X_{t,x}(s^-))-f(X_{t,y}(s^-))
             \rangle
             \big]
       \,\mbox{d}s
      \\&+
      \int_t^{t+\theta}
          \E\Big[
            \int_Z
                |\sigma(X_{t,x}(s^-),z)-\sigma(X_{t,y}(s^-),z)|^2
            \,\nu(\mbox{d}z)
            \Big]
      \mbox{d}s
      \\&-
      2\int_{t}^{t+\theta}
           \E\big[
             \langle
             x-y,f(X_{t,x}(s^-))-f(X_{t,y}(s^-))
             \rangle
             \big]
       \,\mbox{d}s.
\end{split}
\end{equation*}
With the aid of \eqref{eq:gmc} and the Schwarz inequality, one sees that
\begin{equation}\label{eq:txy.ugmc}
\begin{split}
      \E\big[ |\Phi_{t,x,y}(t+\theta)|^2 \big]
      \leq&
      K\int_t^{t+\theta}
           \E\big[ |X_{t,x}(s^-)-X_{t,y}(s^-)|^2 \big]
       \,\mbox{d}s
      \\&+
      2|x-y|
      \int_t^{t+\theta}
          \E\big[ |f(X_{t,x}(s^-))-f(X_{t,y}(s^-)) | \big]
      \,\mbox{d}s.
\end{split}
\end{equation}
It follows from \eqref{eq:anal.solu.bounded} with $q \leq \bar{p}$, \eqref{eq:fpg}, \eqref{eq:tx} and the H\"{o}lder inequality that
\begin{equation}\label{eq:txy.third.ufsg}
\begin{split}
      &
      \int_t^{t+\theta}
          \E\big[ |f(X_{t,x}(s^-))-f(X_{t,y}(s^-))| \big]
      \,\mbox{d}r
      \\\leq&
      \int_t^{t+\theta}
          \E\big[
            \big( 1+|X_{t,x}(s^-)|^q+|X_{t,y}(s^-)|^q \big)^\frac{1}{2}
            \big| X_{t,x}(s^-)-X_{t,y}(s^-) \big|
            \big]
      \,\mbox{d}r
      \\\leq&
      \int_t^{t+\theta}
          \big(
               1
               +\E\big[|X_{t,x}(s^-)|^q\big]
               +\E\big[|X_{t,y}(s^-)|^{q}\big]
          \big)^\frac{1}{2}
          \\&~~~\hspace{1em}\times\big(
               \E\big[ |X_{t,x}(s^-)-X_{t,y}(s^-)|^2 \big]
          \big)^{\frac{1}{2}}
      \,\mbox{d}r
      \\\leq&
      K\theta\big(1+|x|^q+|y|^q\big)^{\frac{1}{2}}|x-y|.
\end{split}
\end{equation}
Plugging \eqref{eq:repre}, \eqref{eq:txy.third.ufsg} into \eqref{eq:txy.ugmc} and using \eqref{eq:elementaryinequality}, the properties of the right-continuous with left limits functions \cite[Section 2.9]{applebaum2009levy} give
\begin{equation*}\label{eq:txy.ufsg}
\begin{split}
      \E\big[  |\Phi_{t,x,y}(t+&\theta)|^{2}  \big]
      \leq
      K\theta\big(1+|x|^{q}+|y|^{q}\big)^{\frac{1}{2}}|x-y|^2
      +
      K\int_{t}^{t+\theta}
       \E\big[  |\Phi_{t,x,y}(r^-)|^{2}  \big]
       \,\mbox{d}r
      \\\leq&
      K\theta\big(1+|x|^{q}+|y|^{q}\big)^{\frac{1}{2}}|x-y|^2
      +
      K\int_{t}^{t+\theta}
      \sup_{t \leq s \leq r}
       \E\big[  |\Phi_{t,x,y}(s)|^{2}  \big]
       \,\mbox{d}r,
\end{split}
\end{equation*}
which immediately shows
\begin{equation*}
      \sup_{t \leq s \leq t+\theta}
      \E\big[  |\Phi_{t,x,y}(s)|^{2}  \big]
      \leq
      K\theta\big(1+|x|^{q}+|y|^{q}\big)^{\frac{1}{2}}|x-y|^2
      +
      K\int_{t}^{t+\theta}\hspace{-0.75em}
      \sup_{t \leq s \leq r}
       \E\big[  |\Phi_{t,x,y}(s)|^{2}  \big]
       \,\mbox{d}r.
\end{equation*}
By the Gronwall inequality, we get
\begin{equation*}
      \sup_{t \leq s \leq t+\theta}
      \E\big[  |\Phi_{t,x,y}(s)|^{2}  \big]
      \leq
      K\theta\big(1+|x|^{q}+|y|^{q}\big)^{\frac{1}{2}}|x-y|^2,
      \quad \forall \theta \in [0,h],
\end{equation*}
and consequently
\begin{equation*}
      \E\big[  |\Phi_{t,x,y}(t+\theta)|^{2}  \big]
      \leq
      K\theta\big(1+|x|^{q}+|y|^{q}\big)^{\frac{1}{2}}|x-y|^2,
      \quad \forall \theta \in [0,h],
\end{equation*}
which gives the desired result \eqref{eq:txy} by taking $\theta = h$.
\end{proof}

Equipped with the above lemma, we are ready to build up the fundamental mean-square convergence theorem for numerical approximations of \eqref{eq:JSDE}.
\begin{theorem}\label{theorem:convergence}
      Let Assumptions \ref{as:assumption1}, \ref{as:assumption2}, \ref{ass:drift-polynomial} be satisfied and let $\bar{p}$ and $q$ come from \eqref{eq:gcc.b} and \eqref{eq:fpg}, respectively, satisfying $q \leq \bar{p}$. Suppose that the one-step approximation $Y_{t,x}(t+h)$ defined by \eqref{eq:one.step.numer} has the following local orders of accuracy, i.e., there exist $h_{0} \in (0,1)$, $K > 0$, and $\alpha \geq 1$ such that for any $h \in (0,h_{0}]$, $t + h \leq T, x \in \R^{d}$,
      \begin{align}
            \label{eq:ed}
            \big| \E\big[X_{t,x}(t+h)-Y_{t,x}(t+h)\big] \big|
            \leq&
            K\big(1+|x|^{2\alpha}\big)^{\frac{1}{2}}h^{p_{1}},
            \\
            \label{eq:msd}
            \big(\E\big[|X_{t,x}(t+h)-Y_{t,x}(t+h)|^2\big]\big)^{\frac{1}{2}}
            \leq&
            K\big(1+|x|^{2\alpha}\big)^{\frac{1}{2}}h^{p_{2}}
      \end{align}
      with
      \begin{equation}\label{eq:parameter}
            p_{2}\geq\tfrac{1}{2},\quad p_1\geq{p_{2}+\tfrac{1}{2}}.
      \end{equation}
      Moreover, the approximation $\{Y_n\}_{0 \leq n \leq N}$ produced by \eqref{eq:one.step.numer} has finite $p$-th moments, i.e., for sufficiently large $p \geq 2$ there exist $ K >0, \beta > 0$ such that for any $h \in (0,h_{0}]$,
      \begin{equation}\label{eq:numerbound}
             \sup_{ 0 \leq n \leq N }
            \E\big[|Y_{n}|^{p}\big]
            \leq
            K\big( 1+\big(\E[|X_0|^{\bar{p}}]\big)^{\beta} \big).
      \end{equation}
      Then there exists $\gamma > 0$ and $K > 0$ independent of $h$ such that
      \begin{equation*}\label{eq:convergence}
            \sup_{ 0 \leq n \leq N }
            \left( \E[|X(t_{n}) - Y_{n}|^2 ] \right)^{\frac{1}{2}}
            \leq
            K\big( 1+(\E[|X_{0}|^{\bar{p}}])^{\gamma} \big) h^{p_{2}-\frac{1}{2}}.
      \end{equation*}
\end{theorem}
\begin{proof} 
Since
\begin{equation*}\label{eq:c.equal}
\begin{split}
      X(t_{n+1}) -& Y_{n+1}
      =
      X_{t_0,X_{0}}(t_{n+1})-Y_{t_0,Y_0}(t_{n+1})
      =
      X_{t_{n},X(t_n)}(t_{n+1})-Y_{t_n,Y_{n}}(t_{n+1})
      \\=&
      \big( X_{t_n,X(t_n)}(t_{n+1})-X_{t_n,Y_n}(t_{n+1}) \big)
      +
      \big( X_{t_n,Y_n}(t_{n+1})-Y_{t_n,Y_n}(t_{n+1}) \big),
\end{split}
\end{equation*}
one can arrive at
\begin{equation}\label{eq:c.mse}
\begin{split}
      \E\big[ &|X(t_{n+1})-Y_{n+1}|^2 \big]
   =
      \E\big[
        |X_{t_n,X(t_n)}(t_{n+1})-X_{t_n,Y_n}(t_{n+1})|^{2}
        \big]
    \\&+
      2\E\big[
         \langle
         X_{t_n,X(t_n)}(t_{n+1})-X_{t_n,Y_n}(t_{n+1})
         ,
         X_{t_n,Y_n}(t_{n+1})-Y_{t_n,Y_n}(t_{n+1})
         \rangle
         \big]
   \\&+
      \E\big[
        |X_{t_n,Y_n}(t_{n+1})-Y_{t_n,Y_n}(t_{n+1})|^2
        \big]
    :=
      A_{1} + A_{2} + A_{3}.
\end{split}
\end{equation}
Thanks to the conditional version of \eqref{eq:tx}, we have
\begin{equation*}\label{eq:c.a1}
      A_{1}
      =
      \E\big[
        \E\big(
               |X_{t_n,X(t_n)}(t_{n+1})-X_{t_n,Y_n}(t_{n+1})|^{2}
               \big|
               \F_{t_n}
          \big)
        \big]
      \leq
      (1+Kh)\E\big[ |X(t_{n})-Y_{n}|^{2} \big].
\end{equation*}
Likewise, the conditional version of \eqref{eq:msd} ensures that
\begin{equation*}\label{eq:c.a2}
      A_{3}
      =
      \E\big[
        \E\big(
                | X_{t_n,Y_n}(t_{n+1})-Y_{t_n,Y_n}(t_{n+1}) |^{2}
                \big|
                \F_{t_{n}}
          \big)
        \big]
      \leq
      Kh^{2p_{2}}
      \E\big[ (1 + |Y_n|^{2\alpha}) \big].
\end{equation*}
It remains to estimate $A_{2}$. Using \eqref{eq:repre} leads to
\begin{equation}\label{eq:c.a3}
\begin{split}
      A_{2}
      =&
      2\E\big[
         \langle
         \Phi_{t_n,X(t_n),Y_n}(t_{n+1})
         ,
         X_{t_n,Y_n}(t_{n+1})-Y_{t_n,Y_n}(t_{n+1})
         \rangle
         \big]
      \\&+
      2\E\big[
         \langle
         X(t_n)-Y_n,X_{t_n,Y_n}(t_{n+1})-Y_{t_n,Y_n}(t_{n+1})
         \rangle
         \big]
      := A_{21} + A_{22}.
\end{split}
\end{equation}
Before estimating $A_{21}$, we put the conditional versions of \eqref{eq:txy} and \eqref{eq:msd} here,
\begin{equation*}\label{eq:312}
      \E\big(|\Phi_{t_{n},X(t_{n}),Y_{n}}(t_{n}+h)|^{2}
      \big| \F_{t_{n}}\big)
      \leq
      Kh\big(1+|X(t_{n})|^q+|Y_{n}|^q\big)^{\frac{1}{2}}
      \big|X(t_{n})-Y_{n}\big|^2,
\end{equation*}
\begin{equation*}\label{eq:323}
      \E\big(|X_{t_{n},Y_{n}}(t_{n}+h)-Y_{t_{n},Y_{n}}(t_{n}+h)|^{2} \big| \F_{t_{n}} \big)
      \leq
      K\big(1+|Y_{n}|^{2\alpha}\big)h^{2p_{2}}.
\end{equation*}
By the Schwarz inequality, the conditional version of the H\"older inequality and the $\F_{t_n}$-measurability of $X(t_n)-Y_n$, we get
\begin{equation*}\label{eq:c.a31}
\begin{split}
      &A_{21}
      =
      2\E\big[
         \E\big(
         \big\langle
             \Phi_{t_n,X(t_n),Y_n}(t_{n+1})
             ,
             X_{t_n,Y_n}(t_{n+1})-Y_{t_n,Y_n}(t_{n+1})
         \big\rangle
           \big|\F_{t_n}\big)
         \big]
      \\\leq&
      2\E\big[
         \big(
              \E(
                      |\Phi_{t_n,X(t_n),Y_n}(t_{n+1})|^2
                      |
                      \F_{t_n}
                )
         \big)^{\frac{1}{2}}
         \big(
              \E\big(
                     |X_{t_n,Y_n}(t_{n+1})-Y_{t_n,Y_n}(t_{n+1})|^2
                     \big|
                     \F_{t_n}
                \big)
         \big)^{\frac{1}{2}}
         \big]
      \\\leq&
      Kh^{p_{2}+\frac{1}{2}}
      \E\big[
          (1+|X(t_{n})|^{q}+|Y_{n}|^{q})^{\frac{1}{4}}
          |X(t_n)-Y_{n}|(1+|Y_{n}|^{2\alpha})^{\frac{1}{2}}
        \big]
      \\=&
      K\E\big[h^{\frac{1}{2}}|X(t_n)-Y_{n}|
         \times
         h^{p_{2}}(1+|X(t_{n})|^{q}+|Y_{n}|^{q})^{\frac{1}{4}}
         (1+|Y_{n}|^{2\alpha})^{\frac{1}{2}}
         \big]
      \\\leq&
      Kh\E\big[ |X(t_n)-Y_n|^{2} \big]
      +
      Kh^{2p_{2}}
      \E\big[
        \big( 1+|Y_n|^{2\alpha} \big)
        \big( 1+|X(t_n)|^q+|Y_n|^q \big)^{\frac{1}{2}}
        \big],
\end{split}
\end{equation*}
where the weighted Young inequality $ab \leq \varepsilon a^{2} + \frac{b^{2}}{4\varepsilon}$ for all $a,b \in \R$ with $\varepsilon = \frac{1}{2}>0$ was used in the last step.
The $\F_{t_n}$-measurability of $X(t_n)-Y_n$, $p_{1} \geq p_{2} + \frac{1}{2}$, the H\"{o}lder inequality and the conditional version of \eqref{eq:ed} imply
\begin{equation*}\label{eq:c.a32.a}
\begin{split}
      A_{22}
      =&
      2\E\big[
      \big\langle
         X(t_n)-Y_n,
         \E\big(
                 X_{t_n,Y_n}(t_{n+1})-Y_{t_n,Y_n}(t_{n+1})
                 \big|
                 \F_{t_n}
           \big)
      \big\rangle
         \big]
      \\\leq&
      2\big(
      \E\big[ |X(t_n)-Y_n|^2 \big]\big)^{\frac{1}{2}}
      \big(
         \E\big[
           \big|
           \E(
           (X_{t_n,Y_n}(t_{n+1})-Y_{t_n,Y_n}(t_{n+1})) | \F_{t_n}
             )
           \big|^{2}
           \big]
      \big)^{\frac{1}{2}}
      \\\leq&
      Kh^{p_1}
      \big(
          \E\big[ |X(t_n)-Y_n|^2 \big]
      \big)^{\frac{1}{2}}
      \E\big[ 1+|Y_n|^{2\alpha} \big]^{\frac{1}{2}}.
      \\\leq&
      Kh^{p_2+\frac{1}{2}}
      \big( \E\big[ |X(t_n)-Y_n |^2 \big] \big)^{\frac{1}{2}}
      \E\big[ 1+|Y_n|^{2\alpha} \big]^{\frac{1}{2}}
      \\\leq&
      Kh\E\big[ |X(t_n)-Y_n|^2 \big]
      +
      Kh^{2p_2} \E\big[ 1+|Y_n|^{2\alpha} \big].
\end{split}
\end{equation*}
Inserting $A_{21}$ and $A_{22}$ into \eqref{eq:c.a3}, we get
\begin{equation*}
      A_{2}
      \leq
      Kh\E\big[ |X(t_n)-Y_n|^2 \big]
      +
      Kh^{2p_{2}}
      \E\big[
        (1+|Y_n|^{2\alpha})
        \big(1+|X(t_n)|^q+|Y_n|^q\big)^{\frac{1}{2}}
        \big].
\end{equation*}
Substituting $A_1, A_2, A_3$ into \eqref{eq:c.mse} and using \eqref{eq:anal.solu.bounded} with $q \leq \bar{p}$, \eqref{eq:numerbound} tell that
\begin{equation*}\label{eq:c.unb}
\begin{split}
      \E\big[ |X(t_{n+1}) - Y_{n+1} |^2\big]
      \leq&
      (1+Kh)\E\big[|X(t_n)-Y_n|^2\big]
      +
      Kh^{2p_2}
      \big(1+\big(\E[|X_{0}|^{\bar{p}}]\big)^{\beta+\frac{q}{2\bar{p}}}\big),
\end{split}
\end{equation*}
which immediately gives
\begin{equation*}
\begin{split}
      &\E\big[ |X(t_{n}) - Y_{n} |^2\big]
      -
      \E\big[|X(t_{n-1})-Y_{n-1}|^2\big]
      \\\leq&
      Kh\E\big[|X(t_{n-1})-Y_{n-1}|^2\big]
      +
      Kh^{2p_2}
      \big(1+\big(\E[|X_{0}|^{\bar{p}}]\big)^{\beta+\frac{q}{2\bar{p}}}\big).
\end{split}
\end{equation*}
By summation, $X(t_{0}) = X_0 = Y_0$ and $nh \leq T$, we have
\begin{equation*}
\begin{split}
      \E\big[ |X(t_{n}) - Y_{n} |^2\big]
      \leq
      Kh\sum\limits_{i=1}^{n-1}\E\big[|X(t_{i})-Y_{i}|^2\big]
      +
      Kh^{2p_2-1}
      \big(1+\big(\E[|X_{0}|^{\bar{p}}]\big)^{\beta+\frac{q}{2\bar{p}}}\big).
\end{split}
\end{equation*}
Exploiting the discrete Gronwall inequality (see, e.g., \cite[Lemma 3.4]{Mao13b}) and using $nh \leq T$ again result in the desired result. 
\end{proof}

\section{Application of the fundamental convergence theorem:
         convergence rates of the tamed Euler method}
\label{sec:tamed}

As the first application of the fundamental mean-square convergence theorem,
we shall construct a new version of the tamed Euler method, also named the tamed Euler method, as follows
\begin{equation}\label{eq:tame.method}
      Y_{n+1}
      =
      Y_n
      +
      \frac{f(Y_n)h}{1+|f(Y_n)|h}
      +
      \frac{g(Y_n)\Delta{W_n}}{1+|g(Y_n)|h}
      +
      \int_{t_n}^{t_{n+1}}\hspace{-0.5em}\int_Z
          \frac{\sigma(Y_n,z)}{1+|\sigma(Y_n,z)|h}
      \,\bar{N}(\mbox{d}s,\mbox{d}z)
\end{equation}
with $Y_0=X_0$, where $\Delta{W_n} := W(t_{n+1})-W(t_n),n=0,1,\ldots,N-1$.
The scheme is different from schemes introduced by Tretyakov and Zhang \cite{Tretyakov13} even the jump term vanishes, i.e., $\sigma \equiv 0$. To apply
Theorem \ref{theorem:convergence}, it is crucial to obtain the boundedness of the higher-order moments of $\{Y_n\}_{0\leq n \leq N}$ given by \eqref{eq:tame.method}.

\subsection{Bounded $p$-th moments of the tamed Euler method}
\label{subset:bounded-moments}

We will present some lemmas before showing the boundedness of $p$-th moments of $\{Y_n\}_{0\leq n \leq N}$.
The first one is 
the Burkholder-Davis-Gundy (BDG) inequality
(see \cite[Lemma 1]{Mikulevicius12}).

\begin{lemma}[Burkholder-Davis-Gundy inequality]\label{lemma:bdg}
      Suppose $p\geq1$ and let $\mathcal{P}$ be the progressive $\sigma$-algebra on $[0,\infty)\times\Omega$
      and $\mathfrak{B}(Z)$ be the Borel $\sigma$-algebra of $Z$. If $\phi $ is a $\mathcal{P}\otimes\mathfrak{B}(Z)$-measurable function such that $\P$-a.s.
      $
          \int_{0}^{T}\int_{Z}
              |\phi(s,z)|^{2}
          \, \nu(\mbox{d}z)\mbox{d}s
          <
          \infty.
      $
Then there exists $K>0$ such that for all $p \geq 2$,
\begin{equation}\label{eq:bdg}
\begin{split}
      \E\bigg[
        \sup_{0\leq{t}\leq{T}}
        \Big|
        \int_{0}^{t}\int_{Z} \phi (s,z) \,\bar{N}(\mbox{d}s,\mbox{d}z)
        \Big|^{p}
        \bigg]
      \leq&
     K\E\Big[
        \int_{0}^{T}\int_{Z}
            |\phi (s, z)|^{2}
        \nu(\mbox{d}z)\mbox{d}s
        \Big]^{\frac{p}{2}}
     \\&+
      K\E\Big[
         \int_{0}^{T}\int_{Z}
             | \phi( s, z ) |^{p}
         \nu(\mbox{d}z)\,\mbox{d}s
         \Big].
\end{split}
\end{equation}
Moreover, if $1 \leq p < 2$, then the last term of \eqref{eq:bdg} can be omitted.
\end{lemma}

Unlike the BDG inequality for the Wiener process, the $p$-th moments ($p \geq 1$) of the Poisson increments
$
  \int_{t}^{t+h} \int_Z \bar{N}(\mbox{d}s,\mbox{d}z), t \geq 0, h > 0
$
contribute to magnitude not more than $O(h)$. This, as already discussed earlier, causes significant difficulties in proving bounded $p$-th moments of the tamed Euler method.  Additionally, we need the following elementary inequality.
\begin{lemma}\label{lemma:varint.binomial}
      Let $1 \leq l \in \N $. Then there exists $K : = K(l) > 0$ such that
      \begin{equation*}
            \big| |x|^{2l}-|y|^{2l} \big|
            \leq
            K\sum\limits_{i=1}^{2l}
             |x-y|^{i}|y|^{{2l}-i},
            \quad
            \forall x,y \in \R^d,
      \end{equation*}
where, as a conventional notation, we set $y^0 = 1$.
\end{lemma}

\begin{proof}
We derive from the binomial formula and \eqref{eq:elementaryinequality} that for all $x,y \in \R^d$,
\begin{equation*}
\begin{split}
      &\big| |x|^{2l}-|y|^{2l} \big|
      =
      \big| \langle x-y+y, x-y+y \rangle^{l} -  |y|^{2l} \big|
      \\=&
      \big|
      \big(
            |x - y|^2
            +
            2\langle x-y, y \rangle + |y|^2
      \big)^l - |y|^{2l}
      \big|
      \\=&
      \Big|
          \sum\limits_{i=1}^{l} \mbox{C}_{l}^{i}
          \big( |x - y|^2 + 2\langle x-y, y \rangle \big)^{i}
          |y|^{2(l-i)}
      \Big|
      \\\leq&
      K\sum\limits_{i=1}^{l}
      \Big(
          |x - y|^{2i} |y|^{2l - 2i}
          +
          |x - y|^{i} |y|^{2l - i}
      \Big)
      \\\leq&
      K\sum\limits_{i=1}^{2l} |x-y|^{i} |y|^{2 l -i},
\end{split}
\end{equation*}
which completes the proof.
\end{proof}

We will prove that the numerical approximations produced by \eqref{eq:tame.method} enjoy bounded high-order moments. At first, we show that the boundedness of high-order moments remains valid within a family of appropriate subevents. Before doing so, we would like to add some comments here.

\begin{remark}\label{remark:Mom-Bound}
      It is worthwhile to emphasize that, one can not simply extend the analysis in \cite{Tretyakov13} to the present jump setting, because the property of Wiener increments $\E [ \| W(t+h)-W(t) \|^l ] = O(h^{\frac{l}{2}}), l \in \N $ was essentially used there (see the treatment of the last term of (3.6) in \cite{Tretyakov13}) while, as clarified earlier, the Poisson increments violate such nice property and the jump coefficients might grow super-linearly.
\end{remark}

To overcome the above difficulty, we work with continuous-time approximations and do very careful estimates. Let $R > 0$ be sufficiently large and define a sequence of decreasing subevents
\begin{equation*}\label{Eq:tame.Omega}
      \Omega_{R,n}
      :=
      \{
        \omega\in\Omega:
        \sup_{ 0 \leq i \leq n }
        |Y_i(\omega)|
        \leq
        R
      \},
      \quad
      \forall n = 0,1,\ldots,N, \,\, N \in \N.
\end{equation*}
Obviously, $\1_{\Omega_{R(h),n}}$ are $\F_{t_n}$-measurable for all $n=0,1,\ldots,N$. The next result indicates that moments of $\{Y_n\}_{0\leq n \leq N}$ are bounded on the  subevents $\Omega_{R,n}$.

\begin{lemma}\label{lemma:tame.moment.bound.a}
      Suppose Assumptions \ref{as:assumption1}, \ref{as:assumption2}, \ref{ass:drift-polynomial} hold. Let $\bar{p} \geq 2$ coming from \eqref{eq:gcc.b} be a sufficiently large even number and let $\{Y_n\}_{0\leq n \leq N}$ be given by \eqref{eq:tame.method}. Then there exist $R = h^{-H(q)^{-1}} := R(h)$ with $H(q)=\max\left\{ 1+q,\tfrac{3}{2}q \right\}$ and $K>0$ independent of $h$ such that
      \begin{alignat}{1}
            \sup_{0 \leq n \leq N}
            \E\big[ \1_{\Omega_{R(h),n}} |Y_{n}|^{\bar{p}} \big]
            &\leq
            K\big(1+\E\big[|X_{0}|^{\bar{p}}\big]\big),
            \label{eq:tame.parta}
            \\
            \sup_{0 \leq n \leq N-1}
            \E\big[ \1_{\Omega_{R(h),n}}|Y_{n+1}|^{\bar{p}} \big]
            &\leq
            K\big(1+\E\big[|X_{0}|^{\bar{p}}\big]\big).
            \label{eq:tame.special.b}
      \end{alignat}
\end{lemma}

\begin{proof}
We define a continuous-time version $\{\bar{Y}(t)\}_{0\leq{t}\leq{T}}$ of $\{Y_n\}_{0\leq{n}\leq{N}}$ by
\begin{equation}\label{eq:tame.contiuous.method.c}
\begin{split}
      \bar{Y}(t)
      =&
      Y_n
      +
      \int_{t_n}^{t}
          \frac{f(Y_n)}{1+|f(Y_n)|h}
      \,\mbox{d}s
      +
      \int_{t_n}^{t}
          \frac{g(Y_n)}{1+|g(Y_n)|h}
      \,\mbox{d}W(s)
      \\&+
      \int_{t_n}^{t}\int_{Z}
          \frac{\sigma(Y_n,z)}{1+|\sigma(Y_n,z)|h}
      \,\bar{N}(\mbox{d}s,\mbox{d}z),
       ~~\P\text{-a.s.}
\end{split}
\end{equation}
for all $t \in [t_n,t_{n+1}], n= 0,1,\ldots,N-1$.
Applying Lemma \ref{lemma:itoformula} yields
\begin{equation*}\label{eq:tame.useIto}
\begin{split}
      |\bar{Y}(t)|^{\bar{p}}
      =&
      |Y_n|^{\bar{p}}
      +
      {\bar{p}}\int_{t_n}^t
           |\bar{Y}(s^-)|^{{\bar{p}}-2}
           \Big\langle
               \bar{Y}(s^-),\frac{f(Y_n)}{1+|f(Y_n)|h}
           \Big\rangle
       \,\mbox{d}s
      \\&+
      \frac{{\bar{p}}}{2}
      \int_{t_n}^{t}
          |\bar{Y}(s^-)|^{{\bar{p}}-2}
          \Big|  \frac{g(Y_n)}{1+|g(Y_n)|h}  \Big|^2
      \,\mbox{d}s
      \\&+
      \frac{{\bar{p}}({\bar{p}}-2)}{2}
      \int_{t_n}^{t}
          |\bar{Y}(s^-)|^{{\bar{p}}-4}
          \Big| \frac{\bar{Y}(s^-)^{T}g(Y_n)}{1+|g(Y_n)|h} \Big|^2
      \,\mbox{d}s
      \\&+
      {\bar{p}}\int_{t_n}^{t}
           |\bar{Y}(s^-)|^{{\bar{p}}-2}
           \Big\langle
               \bar{Y}(s^-),\frac{g(Y_n)\,\mbox{d}W(s)}{1+|g(Y_n)|h}
           \Big\rangle
      \\&+
      \int_{t_n}^{t}\int_{Z}
          \Big|\bar{Y}(s^-)+\frac{\sigma(Y_n,z)}{1+|\sigma(Y_n,z)|h}\Big|^{{\bar{p}}}
          -
          |\bar{Y}(s^-)|^{{\bar{p}}}
      \,N(\mbox{d}s,\mbox{d}z)
      \\&-
      {\bar{p}}\int_{t_n}^{t}\int_{Z}
           |\bar{Y}(s^-)|^{{\bar{p}}-2}
           \Big\langle
               \bar{Y}(s^-),\frac{\sigma(Y_n,z)}{1+|\sigma(Y_n,z)|h}
           \Big\rangle
      \,\nu(\mbox{d}z)\mbox{d}s,
      \quad \P\text{-a.s.}.
\end{split}
\end{equation*}
Then we use the Schwarz inequality
and Lemma \ref{lemma:inequality} to get
\begin{equation*}
\begin{split}
      |\bar{Y}(t)|^{\bar{p}}
      \leq&
      |Y_n|^{\bar{p}}
      +
      {\bar{p}}\int_{t_n}^{t}
           |\bar{Y}(s^-)|^{{\bar{p}}-2}
           \Big\langle
               \bar{Y}(s^-),\frac{f(Y_n)}{1+|f(Y_n)|h}
           \Big\rangle
       \,\mbox{d}s
      \\&+
      \frac{\bar{p}(\bar{p}-1)}{2}
      \int_{t_n}^{t}
          |\bar{Y}(s^-)|^{\bar{p}-2}
          \Big| \frac{g(Y_n)}{1+|g(Y_n)|h} \Big|^{2}
      \,\mbox{d}s
      \end{split}
\end{equation*}\begin{equation*}\label{eq:tame.usecsy}
\begin{split}&+
      \bar{p}\int_{t_n}^{t}
           |\bar{Y}(s^-)|^{\bar{p}-2}
           \Big\langle
               \bar{Y}(s^-),\frac{g(Y_n)\,\mbox{d}W(s)}{1+|g(Y_n)|h}
           \Big\rangle
      +
      K\int_{t_n}^{t} |\bar{Y}(s^-)|^{\bar{p}} \,\mbox{d}s
      \\&+
      \int_{t_n}^{t}\int_{Z}
          \Big|\bar{Y}(s^-)+\frac{\sigma(Y_n,z)}{1+|\sigma(Y_n,z)|h}\Big|^{\bar{p}}
          -
          |\bar{Y}(s^-)|^{\bar{p}}
      \,\bar{N}(\mbox{d}s,\mbox{d}z)
      \\&+
      \big(1+(\bar{p}-2)\varepsilon\big)
      \int_{t_n}^{t}\int_{Z}
         \Big| \frac{\sigma(Y_n,z)}{1+|\sigma(Y_n,z)|h} \Big|^{\bar{p}}
      \,\nu(\mbox{d}z)\mbox{d}s,
      \quad \P\text{-a.s.}.
\end{split}
\end{equation*}
As a result,
\begin{equation*}\label{eq:tame.compensate}
\begin{split}
      &\E\big[ \1_{\Omega_{R,n}}|\bar{Y}(t)|^{\bar{p}} \big]
      \\\leq&
      \E\big[ \1_{\Omega_{R,n}}|Y_n|^{\bar{p}} \big]
      +
      \bar{p}\int_{t_n}^{t}
           \E\Big[
             \1_{\Omega_{R,n}}
             |\bar{Y}(s^-)|^{\bar{p}-2}
             \Big\langle
                 \bar{Y}(s^-),\frac{f(Y_n)}{1+|f(Y_n)|h}
             \Big\rangle
             \Big]
       \,\mbox{d}s
      \\&+
      \frac{\bar{p}(\bar{p}-1)}{2}\hspace{-0.25em}
      \int_{t_n}^t
          \E\big[
            \1_{\Omega_{R,n}}
            |\bar{Y}(s^-)|^{\bar{p}-2}
            |g(Y_n)|^2
            \big]
      \,\mbox{d}s
      +
      K\int_{t_n}^{t}
           \E\big[
             \1_{\Omega_{R,n}} |\bar{Y}(s^-)|^{\bar{p}}
             \big]
       \,\mbox{d}s
      \\&+
      \big(1+(\bar{p}-2)\varepsilon\big)
      \int_{t_n}^{t}
          \E\Big[
            \1_{\Omega_{R,n}}
            \int_{Z}
                |\sigma(Y_n,z)|^{\bar{p}}
            \,\nu(\mbox{d}z)
            \Big]
      \mbox{d}s
      \\=&
      \E\big[\1_{\Omega_{R,n}}|Y_n|^{\bar{p}}\big]
      +
      \bar{p}\int_{t_n}^{t}
           \E\big[
             \1_{\Omega_{R,n}}|Y_n|^{\bar{p}-2}
             \langle Y_n,f(Y_n) \rangle
             \big]
       \,\mbox{d}s
      \\&+
      \frac{\bar{p}(\bar{p}-1)}{2}
      \int_{t_n}^{t}
          \E\big[
            \1_{\Omega_{R,n}} |Y_n|^{\bar{p}-2} |g(Y_n)|^2
            \big]
      \,\mbox{d}s
      +
      K\int_{t_n}^{t}
           \E\big[ \1_{\Omega_{R,n}} |\bar{Y}(s^-)|^{\bar{p}} \big]
       \,\mbox{d}s
      \\&+
      \big(1+(\bar{p}-2)\varepsilon\big)
      \int_{t_n}^{t}
          \E\Big[
            \1_{\Omega_{R,n}}
            \int_{Z} |\sigma(Y_n,z)|^{\bar{p}} \,\nu(\mbox{d}z)
            \Big]
      \mbox{d}s
      +I_1+I_2+I_3+I_4,
\end{split}
\end{equation*}
where
\begin{align*}
      I_{1}:=&
      \bar{p}\int_{t_n}^{t}
           \E\Big[
             \1_{\Omega_{R,n}}\big|\bar{Y}(s^-)\big|^{\bar{p}-2}
             \Big\langle
             \bar{Y}(s^-)-Y_n,\frac{f(Y_n)}{1+|f(Y_n)|h}
             \Big\rangle
             \Big]
       \,\mbox{d}s,
      \\
      I_{2}:=&
      \bar{p}\int_{t_n}^{t}
           \E\Big[
             \1_{\Omega_{R,n}}
             |\bar{Y}(s^-)|^{\bar{p}-2}
             \Big\langle
                 Y_n,\frac{f(Y_n)}{1+|f(Y_n)|h}-f(Y_n)
             \Big\rangle
             \Big]
       \,\mbox{d}s,
      \\
      I_{3}:=&
      \bar{p}\int_{t_n}^{t}
           \E\big[
             \1_{\Omega_{R,n}}
             \big(|\bar{Y}(s^-)|^{\bar{p}-2}-|Y_n|^{\bar{p}-2}\big)
             \langle Y_n,f(Y_n)\rangle
             \big]
       \,\mbox{d}s,
      \\
      I_{4}:=&
      \frac{\bar{p}(\bar{p}-1)}{2}
      \int_{t_n}^{t}
          \E\big[
            \1_{\Omega_{R,n}}
            \big(|\bar{Y}(s^-)|^{\bar{p}-2}-|Y_n|^{\bar{p}-2}\big) |g(Y_n)|^{2}
            \big]
      \,\mbox{d}s.
\end{align*}
Using the coercivity condition \eqref{eq:gcc.b} yields
\begin{equation}\label{eq:tame.using_gcc}
\begin{split}
      \E[ \1_{\Omega_{R,n}} |\bar{Y}(t)|^{\bar{p}} ]
      \leq&
      (1+Kh)\E[ \1_{\Omega_{R,n}} |Y_{n}|^{\bar{p}} ]
      +
      K\int_{t_n}^{t}
           \E[ \1_{\Omega_{R,n}} |\bar{Y}(s^-)|^{\bar{p}} ]
       \,\mbox{d}s
      \\&+I_1+I_2+I_3+I_4.
\end{split}
\end{equation}
For $I_1$, using the Schwarz inequality and \eqref{eq:elementaryinequality} helps us to get
\begin{equation}\label{eq:I.one.1}
\begin{split}
      I_1
      \leq&
      \bar{p}\int_{t_n}^{t}
           \E\big[
             \1_{\Omega_{R,n}}
             |\bar{Y}(s^-)|^{\bar{p}-2}
             |\bar{Y}(s^-)-Y_{n}|
             |f(Y_n)|
             \big]
           \mbox{d}s
      \\\leq&K \!
               \int_{t_n}^t
               \!
                   \E\big[
                       \1_{\Omega_{R,n}}
                       |\bar{Y}(s^-)-Y_n|^{\bar{p}-1}
                       |f(Y_n)|
                     \big]
               \mbox{d}s
      \\&
      +
            K \!
               \int_{t_n}^{t}
               \!
                   \E\big[
                       \1_{\Omega_{R,n}}
                       |\bar{Y}(s^-)-Y_n|
                       |f(Y_n)|
                       |Y_n|^{\bar{p}-2}
                     \big]
               \mbox{d}s
            :=
            I_{11} + I_{12}.
\end{split}
\end{equation}
It follows from \eqref{eq:tame.contiuous.method.c} and \eqref{eq:elementaryinequality}
that
\begin{equation}\label{eq:I.one.one.1}
\begin{split}
      I_{11}
      \leq&
      K\int_{t_n}^{t}
           \E\Big[
             \1_{\Omega_{R,n}} |f(Y_n)|
             \Big| \frac{(s-t_n)f(Y_n)}{1+|f(Y_n)|h} \Big|^{\bar{p}-1}
             \Big]
       \,\mbox{d}s
      \\&+
      K\int_{t_n}^{t}
       \E\Big[
           \1_{\Omega_{R,n}} |f(Y_n)|
           \Big| \frac{g(Y_n)(W(s)-W(t_{n}))}{1+|g(Y_n)|h} \Big|^{\bar{p}-1}
         \Big]
       \,\mbox{d}s
      \\&+
      K\int_{t_n}^{t}
           \E\Big[
             \1_{\Omega_{R,n}} |f(Y_n)|
             \Big|
                \int_{t_n}^s\int_Z
                   \frac{\sigma(Y_n,z)}{1+|\sigma(Y_n,z)|h}
                \,\bar{N}(\mbox{d}r,\mbox{d}z)
             \Big|^{\bar{p}-1}
             \Big]
       \,\mbox{d}s
      \\\leq&
      Kh^{\bar{p}} \E\big[ \1_{\Omega_{R,n}} |f(Y_n)|^{\bar{p}} \big]
      +
      Kh^{\frac{\bar{p}+1}{2}}
      \E\big[
        \1_{\Omega_{R,n}} |f(Y_n)| |g(Y_n)|^{\bar{p}-1}
        \big]
      \\&+
      Kh\E\Big[
        \sup_{t_n\leq{s}\leq{t}}
        \Big|
           \int_{t_n}^s\int_{Z}
               \frac{\1_{\Omega_{R,n}}\sigma(Y_n,z)
               |f(Y_n)|^{\frac{1}{\bar{p}-1}}}{1+|\sigma(Y_n,z)|h}
           \,\bar{N}(\mbox{d}r,\mbox{d}z)
        \Big|^{\bar{p}-1}
        \Big].
\end{split}
\end{equation}
We then use Lemma \ref{lemma:bdg}, \eqref{eq:fsg}, \eqref{eq:gsg} and \eqref{eq:hsg} to obtain
\begin{equation}\label{eq:I.one.one.2}
\begin{split}
      I_{11}
      \leq&
      Kh^{\bar{p}} \E\big[ \1_{\Omega_{R,n}} |f(Y_n)|^{\bar{p}} \big]
      +
      Kh^{\frac{\bar{p}+1}{2}}
      \E\big[
        \1_{\Omega_{R,n}} |f(Y_n)| |g(Y_n)|^{\bar{p}-1}
        \big]
      \\&+
      Kh^{\frac{\bar{p}+1}{2}}
      \E\Big[
        \1_{\Omega_{R,n}} |f(Y_n)|
        \Big(
           \int_Z |\sigma(Y_n,z)|^{2}\,\nu(\mbox{d}z)
        \Big)^{\frac{\bar{p}-1}{2}}
        \Big]
      \\&+
      Kh^{2}
      \E\Big[
        \1_{\Omega_{R,n}} |f(Y_n)|
        \int_Z |\sigma(Y_n,z)|^{\bar{p}-1} \,\nu(\mbox{d}z)
        \Big]
      \\\leq&
      Kh
      +
      Kh^{\bar{p}}
      \E\big[ \1_{\Omega_{R,n}} |Y_n|^{\bar{p}+\frac{\bar{p}q}{2}} \big]
      +
      Kh^{\frac{\bar{p}+1}{2}}
      \E\big[ \1_{\Omega_{R,n}} |Y_n|^{\bar{p}+\frac{(\bar{p}+1)q}{4}} \big]
      \\&+
      Kh^{2} \E\big[ \1_{\Omega_{R,n}} |Y_n|^{\bar{p}+q} \big].
\end{split}
\end{equation}
Repeating the same arguments used in \eqref{eq:I.one.one.1} and \eqref{eq:I.one.one.2} implies
\begin{equation*}
\begin{split}
      I_{12}
      \leq&
      \int_{t_n}^{t}
          \E\Big[
            \1_{\Omega_{R,n}} |f(Y_n)||Y_n|^{\bar{p}-2}
            \Big| \frac{(s-t_n)f(Y_n)}{1+|f(Y_n)|h} \Big|
            \Big]\,\mbox{d}s
      \\&+
      \int_{t_n}^{t}
          \E\Big[
              \1_{\Omega_{R,n}} |f(Y_n)||Y_n|^{\bar{p}-2}
              \Big| \frac{g(Y_n)(W(s)-W(t_{n}))}{1+|g(Y_n)|h} \Big|
            \Big]
      \,\mbox{d}s
      \\&+
      \int_{t_n}^{t}
          \E\Big[
            \1_{\Omega_{R,n}} |f(Y_n)| |Y_n|^{\bar{p}-2}
            \Big|
               \int_{t_n}^s\int_Z
                  \frac{\sigma(Y_n,z)}{1+|\sigma(Y_n,z)|h}
               \,\bar{N}(\mbox{d}r,\mbox{d}z)
            \Big|
            \Big]
      \,\mbox{d}s
\end{split}
\end{equation*}
and thus
\begin{equation*}\label{eq:I.one.two}
\begin{split}
      I_{12}
      \leq&
      h^{2}\E\big[ \1_{\Omega_{R,n}} |f(Y_n)|^{2} |Y_n|^{\bar{p}-2} \big]
      +
      h^{\frac{3}{2}}
      \E\big[
        \1_{\Omega_{R,n}} |g(Y_n)| |f(Y_n)| |Y_n|^{\bar{p}-2}
        \big]
      \\&+
      h\E\Big[
         \sup_{t_n\leq{s}\leq{t}}
         \Big|
            \int_{t_n}^s\int_Z
                \frac{\1_{\Omega_{R,n}}\sigma(Y_n,z)|f(Y_n)||Y_n|^{\bar{p}-2}}
                     {1+|\sigma(Y_n,z)|h}
            \,\bar{N}(\mbox{d}r,\mbox{d}z)
         \Big|
         \Big]
      \\\leq&
      h^{2}
      \E\big[
        \1_{\Omega_{R,n}} |f(Y_n)|^{2} |Y_n|^{\bar{p}-2}
        \big]
      +
      h^{\frac{3}{2}}
      \E\big[
        \1_{\Omega_{R,n}} |g(Y_n)| |f(Y_n)| |Y_n|^{\bar{p}-2}
        \big]
      \\&+
      Kh^{\frac{3}{2}}
      \E\Big[
        \1_{\Omega_{R,n}} |f(Y_n)| |Y_n|^{\bar{p}-2}
        \Big(\int_Z |\sigma(Y_n,z)| \,\nu(\mbox{d}z)\Big)^{\frac{1}{2}}
        \Big]
      \\\leq&
      Kh
      +
      Kh^{2}\E\big[ \1_{\Omega_{R,n}}|Y_n|^{\bar{p}+q} \big]
      +
      Kh^{\frac{3}{2}}
      \E\big[ \1_{\Omega_{R,n}}|Y_n|^{\bar{p}+\frac{3q}{4}} \big].
\end{split}
\end{equation*}
Substituting the above $I_{11}$ and $I_{12}$ into \eqref{eq:I.one.1} gives
\begin{equation}\label{eq:I.one}
\begin{split}
      I_1
      \leq&
      Kh
      +
      Kh^{\bar{p}}
      \E\big[
        \1_{\Omega_{R,n}} |Y_n|^{\bar{p}+\frac{\bar{p}q}{2}}
        \big]
      +
      Kh^{\frac{\bar{p}+1}{2}}
      \E\big[ \1_{\Omega_{R,n}} |Y_n|^{\bar{p}+\frac{(\bar{p}+1)q}{4}} \big]
      \\&+
      Kh^{2} \E\big[ \1_{\Omega_{R,n}}|Y_n|^{\bar{p}+q} \big]
      +
      Kh^{\frac{3}{2}}
      \E\big[ \1_{\Omega_{R,n}}|Y_n|^{\bar{p}+\frac{3q}{4}} \big].
\end{split}
\end{equation}
Treating $I_2$ by the Schwarz inequality, the Young inequality and \eqref{eq:fsg} leads to
\begin{equation}\label{eq:I.two}
\begin{split}
      I_2
      \leq&
      \bar{p} \int_{t_n}^{t}
           \E\big[
             \1_{\Omega_{R,n}}
             |\bar{Y}(s^-)|^{\bar{p}-2} |Y_n| |f(Y_n)|^{2} h
             \big]
       \,\mbox{d}s
      \\\leq&
      (\bar{p}-2)
      \int_{t_n}^{t}
          \E\big[ \1_{\Omega_{R,n}} |\bar{Y}(s^-)|^{\bar{p}} \big]
      \,\mbox{d}s
      +
      2h^{1+\frac{\bar{p}}{2}}
      \E\big[
        \1_{\Omega_{R,n}} |Y_{n}|^{\frac{\bar{p}}{2}} |f(Y_n)|^{\bar{p}}
        \big]
      \\\leq&
      (\bar{p}-2)
      \int_{t_n}^{t}
          \E\big[ \1_{\Omega_{R,n}} |\bar{Y}(s^-)|^{\bar{p}} \big]
      \,\mbox{d}s
      +
      Kh
      +
      Kh^{1+\frac{\bar{p}}{2}}
      \E\big[ \1_{\Omega_{R,n}} |Y_n|^{\frac{3+q}{2}\bar{p}} \big].
\end{split}
\end{equation}
By the Schwarz inequality, Lemma \ref{lemma:varint.binomial} and \eqref{eq:tame.contiuous.method.c}, we have
\begin{equation*}\label{eq:I.three}
\begin{split}
      &I_3
      \leq
      \bar{p}\int_{t_n}^{t}
           \E\big[
             \1_{\Omega_{R,n}}
             \big| |\bar{Y}(s^-)|^{\bar{p}-2}-|Y_n|^{\bar{p}-2} \big|
             |Y_n| |f(Y_n)|
             \big]
       \,\mbox{d}s
      \\\leq&
      K\sum\limits_{i=1}^{\bar{p}-2}
       \int_{t_n}^{t}
           \E\big[
             \1_{\Omega_{R,n}}
             |\bar{Y}(s^-)-Y_n|^{i} |Y_n|^{\bar{p}-1-i}|f(Y_n)|
             \big]
       \,\mbox{d}s
      \\\leq&
      K\sum\limits_{i=1}^{\bar{p}-2}h^{i+1}
      \E\big[
        \1_{\Omega_{R,n}}
        |f(Y_n)|^{i+1}|Y_n|^{\bar{p}-1-i}
        \big]
      \\&+
      K\sum\limits_{i=1}^{\bar{p}-2}h^{\frac{i}{2}+1}
        \E\big[
          \1_{\Omega_{R,n}}
          |g(Y_n)|^{i}|Y_n|^{\bar{p}-1-i}|f(Y_n)|
          \big]
      \\&+
      Kh\sum\limits_{i=1}^{\bar{p}-2}
        \E\Big[
          \sup_{t_n\leq{s}\leq{t}}
          \Big|
             \int_{t_n}^{s}\hspace{-0.25em}\int_{Z}\hspace{-0.25em}
                 \frac{
                       \1_{\Omega_{R,n}}
                       \sigma(Y_n,z)
                       (|Y_n|^{\bar{p}-1-i}|f(Y_n)|)^{\frac{1}{i}}
                      }
                      {1+|\sigma(Y_n,z)|h}
             \,\bar{N}(\mbox{d}r,\mbox{d}z)
          \Big|^i
          \Big].
\end{split}
\end{equation*}
The techniques used in \eqref{eq:I.one.one.1}--\eqref{eq:I.one.one.2} further tell us that
\begin{equation*}
\begin{split}
      I_3
      \leq&
      K
      \sum\limits_{i=1}^{\bar{p}-2} h^{i+1}
      \E\big[ \1_{\Omega_{R,n}} |f(Y_n)|^{i+1} |Y_n|^{\bar{p}-1-i} \big]
      \\&+
      K\sum\limits_{i=1}^{\bar{p}-2} h^{\frac{i}{2}+1}
      \E\big[
        \1_{\Omega_{R,n}} |g(Y_n)|^i |Y_n|^{\bar{p}-1-i} |f(Y_n)|
        \big]
      \\&+
      K\sum\limits_{i=1}^{\bar{p}-2} h^{\frac{i}{2}+1}
      \E\Big[
        \1_{\Omega_{R,n}}
        |Y_n|^{\bar{p}-1-i}|f(Y_n)|
        \Big(
           \int_Z
               |\sigma(Y_n,z)|^{2}
           \,\nu(\mbox{d}z)
        \Big)^\frac{i}{2}
        \Big]
      \\&+
      Kh^{2}
      \sum\limits_{i=1}^{\bar{p}-2}
      \E\Big[
        \1_{\Omega_{R,n}}
        |Y_n|^{\bar{p}-1-i} |f(Y_n)|
        \int_{Z} |\sigma(Y_n,z)|^{i} \,\nu(\mbox{d}z)
        \Big]
      \\\leq&
      Kh
      +
      K\sum\limits_{i=1}^{\bar{p}-2}h^{i+1}
      \E\big[ \1_{\Omega_{R,n}}|Y_n|^{\bar{p}+\frac{i+1}{2}q} \big]
      +
      K\sum\limits_{i=1}^{\bar{p}-2}h^{\frac{i}{2}+1}
      \E\big[ \1_{\Omega_{R,n}} |Y_n|^{\bar{p}+\frac{i+2}{4}q} \big].
\end{split}
\end{equation*}
Similarly,
\begin{equation*}\label{eq:I.four}
\begin{split}
      I_4
      \leq&
      K\sum\limits_{i=1}^{\bar{p}-2}
       \int_{t_n}^{t}
          \E\big[
            \1_{\Omega_{R,n}}
            |\bar{Y}(s^-)-Y_n|^{i}|g(Y_n)|^2|Y_n|^{\bar{p}-2-i}
            \big]\,\mbox{d}s
      \\\leq&
      Kh
      +
      K\sum\limits_{i=1}^{\bar{p}-2}h^{i+1}
      \E\big[ \1_{\Omega_{R,n}}|Y_n|^{\bar{p}+\frac{1+i}{2}q} \big]
      +
      K\sum\limits_{i=1}^{\bar{p}-2}h^{\frac{i}{2}+1}
      \E\big[ \1_{\Omega_{R,n}}|Y_n|^{\bar{p}+\frac{2+i}{4}q} \big].
\end{split}
\end{equation*}
Inserting $I_{1}$, $I_{2}$, $I_{3}$ and $I_{4}$ into \eqref{eq:tame.using_gcc} promises
\begin{equation*}
\begin{split}
      \E\big[\1_{\Omega_{R,n}}&|\bar{Y}(t)|^{\bar{p}}\big]
      \leq
      Kh
      +
      (1+Kh)
      \E\big[ \1_{\Omega_{R,n}} |Y_n|^{\bar{p}} \big]
      \\&+
      K\int_{t_n}^{t}
           \E\big[ \1_{\Omega_{R,n}} |\bar{Y}(s^-)|^{\bar{p}} \big]
       \,\mbox{d}s
      +
      Kh^{1+\frac{\bar{p}}{2}}
      \E\big[ \1_{\Omega_{R,n}} |Y_n|^{\frac{3+q}{2}\bar{p}} \big]
      \\&+
      K\sum\limits_{i=1}^{\bar{p}-1}h^{i+1}
      \E\big[ \1_{\Omega_{R,n}}|Y_n|^{\bar{p}+\frac{i+1}{2}q} \big]
      +
      K\sum\limits_{i=1}^{\bar{p}-1}h^{\frac{i}{2}+1}
      \E\big[ \1_{\Omega_{R,n}} |Y_n|^{\bar{p}+\frac{i+2}{4}q} \big].
\end{split}
\end{equation*}
Choosing $R = R(h) = h^{-H(q)^{-1}}$ with $H(q)=\max\left\{ 1+q,\tfrac{3}{2}q \right\}$,
one can easily get $\1_{\Omega_{R(h),n}}|Y_{n}| \leq h^{-H(q)^{-1}}$ and for all $i=1,\ldots,\bar{p}-1$,
\begin{alignat*}{1}
      \1_{\Omega_{R(h),n}}|Y_{n}|^{1+q}h
      \leq h^{1-(1+q)H(q)^{-1}} \leq h_{0}^{1-(1+q)H(q)^{-1}},
      \\
      \1_{\Omega_{R(h),n}}|Y_n|^{\frac{1+1/i}{2}q}h
      \leq h^{1-\frac{1+1/i}{2}qH(q)^{-1}} \leq h_{0}^{1-\frac{1+1/i}{2}qH(q)^{-1}},
      \\
      \1_{\Omega_{R(h),n}}|Y_n|^{\frac{1+2/i}{2}q}h
      \leq h^{1-\frac{1+2/i}{2}qH(q)^{-1}} \leq h_{0}^{1-\frac{1+2/i}{2}qH(q)^{-1}}.
\end{alignat*}
These inequalities immediately show that for all $i=1,\ldots,\bar{p}-1$,
\begin{alignat*}{1}
      \1_{\Omega_{R,n}}|Y_n|^{\frac{3+q}{2}\bar{p}}h^{\frac{\bar{p}}{2}}
      =
      \1_{\Omega_{R(h),n}}|Y_n|^{\bar{p}}
      \big(\1_{\Omega_{R(h),n}}|Y_n|^{1+q}h\big)^{\frac{\bar{p}}{2}}
      &\leq
      K\1_{\Omega_{R(h),n}}|Y_n|^{\bar{p}},
      \\
      \1_{\Omega_{R,n}}|Y_n|^{\bar{p}+\frac{i+1}{2}q}h^{i}
      =
      \1_{\Omega_{R(h),n}}|Y_n|^{\bar{p}}
      \big(\1_{\Omega_{R(h),n}}|Y_n|^{\frac{1+1/i}{2}q}h\big)^{i}
      &\leq
      K\1_{\Omega_{R(h),n}}|Y_n|^{\bar{p}},
      \\
      \1_{\Omega_{R,n}}|Y_n|^{\bar{p}+\frac{i+2}{4}q}h^{\frac{i}{2}}
      =
      \1_{\Omega_{R(h),n}}|Y_n|^{\bar{p}}
      \big(\1_{\Omega_{R(h),n}}|Y_n|^{\frac{1+2/i}{2}q}h\big)^{\frac{i}{2}}
      &\leq
      K\1_{\Omega_{R(h),n}}|Y_n|^{\bar{p}},
\end{alignat*}
where the constants $K$ are independent of stepsize $h$. Thus
\begin{equation*}
      \E\big[ \1_{\Omega_{R(h),n}} |\bar{Y}(t)|^{\bar{p}} \big]
      \leq
      Kh
      +
      (1+Kh)
      \E\big[ \1_{\Omega_{R(h),n}} |Y_n|^{\bar{p}} \big]
      +
      K\hspace{-0.25em}\int_{t_n}^{t}
           \E\big[ \1_{\Omega_{R(h),n}} |\bar{Y}(s^-)|^{\bar{p}} \big]
       \,\mbox{d}s.
\end{equation*}
Following the techniques used in \eqref{eq:analbound.gcc}--\eqref{eq:analbound.gcccc}, we have
\begin{equation*}
\begin{split}
      \sup\limits_{t_n \leq r \leq t}
      \E[ \1_{\Omega_{R(h),n}} |\bar{Y}(r)|^{\bar{p}} ]
      \leq&
      Kh
      +
      (1+Kh)
      \E[ \1_{\Omega_{R(h),n}} |Y_n|^{\bar{p}} ]
      \\&+
      K\int_{t_n}^{t}
           \sup\limits_{t_n \leq r \leq s}
           \E[ \1_{\Omega_{R(h),n}} |\bar{Y}(r)|^{\bar{p}} ]
       \,\mbox{d}s.
\end{split}
\end{equation*}
The Gronwall inequality shows that for all $t \in [t_n,t_{n+1}]$,
\begin{equation*}\label{eq:tame.use.cont.gron}
      \sup\limits_{t_n \leq r \leq t}
      \E\big[ \1_{\Omega_{R(h),n}} |\bar{Y}(r)|^{\bar{p}} \big]
      \leq
      Khe^{Kh}
      +
      (1+Kh)e^{Kh}
      \E\big[ \1_{\Omega_{R(h),n}}|Y_n|^{\bar{p}} \big].
\end{equation*}
Taking $t = t_{n+1}$ and repeating the treatment used in \eqref{eq:tx.ugronw} particularly yield
\begin{equation}\label{eq:tame.special}
      \E\big[ \1_{\Omega_{R(h),n}}|Y_{n+1}|^{\bar{p}} \big]
      \leq
      Kh
      +
      (1+Kh)\E\big[ \1_{\Omega_{R(h),n}}|Y_n|^{\bar{p}} \big].
\end{equation}
By the decreasing property of $\{\Omega_{R,n}\}_{0 \leq n \leq N}$, we get $\1_{\Omega_{R(h),n+1}} \leq \1_{\Omega_{R(h),n}}$ and
\begin{equation*}\label{eq:tame.special.a}
      \E\big[ \1_{\Omega_{R(h),n+1}}|Y_{n+1}|^{\bar{p}} \big]
      \leq
      \E\big[ \1_{\Omega_{R(h),n}}|Y_{n+1}|^{\bar{p}} \big]
      \leq
      Kh
      +
      (1+Kh)\E\big[ \1_{\Omega_{R(h),n}}|Y_n|^{\bar{p}} \big],
\end{equation*}
which obviously shows
\begin{equation*}
      \E\big[ \1_{\Omega_{R(h),n}}|Y_{n}|^{\bar{p}} \big]
      \leq
      Knh
      +
      \E\big[ \1_{\Omega_{R(h),0}}|Y_{0}|^{\bar{p}} \big]
      +
      Kh\sum\limits_{i=0}^{n-1}\E\big[ \1_{\Omega_{R(h),i}}|Y_{i}|^{\bar{p}} \big].
\end{equation*}
Applying the discrete Gronwall inequality (see, e.g., \cite[Lemma 3.4]{Mao13b}) and using $nh \leq T$ guarantee \eqref{eq:tame.parta}, which together with  \eqref{eq:tame.special} and $0 < h \leq h_{0}$ immediately suggests \eqref{eq:tame.special.b}. Thus we complete the proof.
\end{proof}

Equipped with the previous lemma, one can derive bounded moments of \eqref{eq:tame.method}.
\begin{lemma}\label{lemma:tame.moment.bound}
      Suppose Assumptions \ref{as:assumption1}, \ref{as:assumption2}, \ref{ass:drift-polynomial} hold and let $\{Y_n\}_{0\leq n \leq N}$ be given by \eqref{eq:tame.method}. Let $H(q)=\max\left\{ 1+q,\tfrac{3}{2}q \right\}$ and let $\bar{p} \geq 2+4H(q)$ be a sufficiently large even number. Then there exist $\beta > 0$ and $K>0$ independent of $h$ such that
  \begin{equation}\label{eq:tame.monent.bound}
      \sup\limits_{ 0 \leq n \leq N }
      \E\big[|Y_n|^p\big]
      \leq
      K\big( 1+\big(\E[|X_0|^{\bar{p}}]\big)^{\beta} \big),
      \quad
      \forall p \in \big[2,\tfrac{\bar{p}-H(q)}{1+\frac{3}{2}H(q)}\big].
  \end{equation}
\end{lemma}


\begin{proof}
It follows from \eqref{eq:tame.method} that
\begin{equation}\label{eq:tame.bound}
\begin{split}
      |Y_{n+1}|&
      \leq
      |Y_n|
      +
      \frac{|f(Y_n)|h}{1+|f(Y_n)|h}
      +
      \frac{|g(Y_n)\Delta{W_n}|}{1+|g(Y_n)|h}
      \\&+
      \Big|
         \int_{t_n}^{t_{n+1}}\int_Z
             \frac{\sigma(Y_n,z)}{1+|\sigma(Y_n,z)|h}
         \,\bar{N}(\mbox{d}s,\mbox{d}z)
      \Big|
      \\\leq&
      |Y_n|
      +
      1
      +
      \frac{|\Delta{W_n}|}{h}
      +
      \nu(Z)
      +
      \int_{t_n}^{t_{n+1}}\int_Z
          \frac{1}{h}
      \,N(\mbox{d}s,\mbox{d}z)
      \leq
      \cdots
      \\\leq&
      |X_0|
      +
      (n+1)(1+\nu(Z))
      +
      h^{-1}\sum\limits_{i=0}^{n}|\Delta{W_i}|
      +
      h^{-1}\int_{0}^{t_{n+1}}\int_Z \,N(\mbox{d}s,\mbox{d}z).
\end{split}
\end{equation}
In view of \eqref{eq:tame.parta}, it suffices to verify $\E\big[\1_{\Omega_{R(h),n}^c}|Y_n|^p\big] < \infty$. Note that
\begin{equation*}\label{eq:rhc}
\begin{split}
       \1_{\Omega_{R(h),n}^{c}}
       &=
       1-\1_{\Omega_{R(h),n}}
       =
       1-\1_{\Omega_{R(h),n-1}}\1_{|Y_n|\leq{R(h)}}
       \\&=
       \1_{\Omega_{R(h),n-1}^{c}}+\1_{\Omega_{R(h),n-1}}\1_{|Y_n|>{R(h)}}
       =
       \sum\limits_{i=0}^{n}\1_{\Omega_{R(h),i-1}}\1_{|Y_{i}|>{R(h)}},
\end{split}
\end{equation*}
where we set $\1_{\Omega_{R(h),-1}}=1$. This together with the H\"{o}lder inequality with $\frac{1}{p'}+\frac{1}{q'}=1$ for
$q'=\frac{\bar{p}}{(\frac{3p}{2}+1)H(q)} > 1$,
due to
$p \leq \tfrac{\bar{p}-H(q)}{1+\frac{3}{2}H(q)}
   < \tfrac{2}{3}(\tfrac{\bar{p}}{H(q)}-1)$,
and the Chebyshev inequality gives
\begin{equation}\label{eq:usesineb}
\begin{split}
      \E\big[\1_{\Omega_{R(h),n}^c}|Y_n|^p\big]
      =&
      \sum\limits_{i=0}^n
      \E\big[
        |Y_n|^p\1_{\Omega_{R(h),i-1}}\1_{|Y_i|>{R(h)}}
        \big]
      \\\leq&
      \sum\limits_{i=0}^{n}
      \big( \E\big[|Y_n|^{pp'}\big] \big)^{\frac{1}{p'}}
      \big(
            \E\big[\1_{\Omega_{R(h),i-1}}\1_{|Y_i|>{R(h)}}\big]
      \big)^{\frac{1}{q'}}
      \\=&
      \big( \E[|Y_n|^{pp'}] \big)^{\frac{1}{p'}}
      \sum\limits_{i=0}^{n}
      \big(
           \P(\1_{\Omega_{R(h),i-1}}|Y_i|>{R(h)})
      \big)^{\frac{1}{q'}}
      \\\leq&
      \big( \E[|Y_n|^{pp'}] \big)^{\frac{1}{p'}}
      \sum\limits_{i=0}^{n}
      \frac{\big(
            \E[\1_{\Omega_{R(h),i-1}}|Y_{i}|^{\bar{p}}]
            \big)^{\frac{1}{q'}}}
           {{R(h)}^{(\frac{3p}{2}+1)H(q)}}.
\end{split}
\end{equation}
Since $p \leq \tfrac{\bar{p}-H(q)}{1+\frac{3}{2}H(q)}$ implies $pp' \leq \bar{p}$, using H\"{o}lder's inequality, \eqref{eq:tame.bound}, \eqref{eq:elementaryinequality} and \eqref{eq:bdg} implies
\begin{equation*}
\begin{split}
      \big( \E[&|Y_n|^{pp'}] \big)^{\frac{1}{p'}}
      \leq
      \big( \E[|Y_n|^{\bar{p}}] \big)^{\frac{p}{\bar{p}}}
      \leq
      4^{p}
      \Big(
              \E\big[|X_0|^{\bar{p}}\big]
              +
              n^{\bar{p}}(1+\nu(Z))^{\bar{p}}
              \\&+
              h^{-\bar{p}}
              \E\big[
                \big|\sum\limits_{i=0}^{n-1} |\Delta{W_i}|\big|^{\bar{p}}
                \big]
              +
              h^{-\bar{p}}
              \E\big[
                \big|\int_{0}^{t_{n}}\int_Z\,N(\mbox{d}s,\mbox{d}z)\big|^{\bar{p}}
                \big]
      \Big)^{\frac{p}{\bar{p}}}
      \\\leq&
      4^{p}
      \Big(
              \E\big[|X_0|^{\bar{p}}\big]
              +
              n^{\bar{p}}(1+\nu(Z))^{\bar{p}}
              +
              h^{-\bar{p}}n^{\bar{p}-1}
              \sum\limits_{i=0}^{n-1}
              \E\big[
                |\Delta{W_i}|^{\bar{p}}
                \big]
              \end{split}
\end{equation*}\begin{equation}\label{eq:usesinecpre1}
\begin{split}&+
              2^{\bar{p}}h^{-\bar{p}}
              \E\big[
                \big|\int_{0}^{t_{n}}\int_Z\,\nu(\mbox{d}z)\mbox{d}s\big|^{\bar{p}}
                \big]
              +
              2^{\bar{p}}h^{-\bar{p}}
              \E\big[
                \big|\int_{0}^{t_{n}}\int_Z\,\bar{N}(\mbox{d}s,\mbox{d}z)\big|^{\bar{p}}
                \big]
      \Big)^{\frac{p}{\bar{p}}}
      \\\leq&
      Kh^{-\frac{3p}{2}}
      +
      K\big( \E[|X_0|^{\bar{p}}] \big)^{\frac{p}{\bar{p}}}.
\end{split}
\end{equation}
Inserting \eqref{eq:usesinecpre1} into \eqref{eq:usesineb} and exploiting the H\"{o}lder inequality, $R(h) = h^{-H(q)^{-1}}$, \eqref{eq:tame.special.b} and \eqref{eq:elementaryinequality}, we deduce
\begin{equation*}
\begin{split}
      \E\big[ \1_{\Omega_{R(h),n}^c}|Y_n|^{p} \big]
      \leq&
      K(n+1)h^{\frac{3p}{2}+1}
      \big(
           Kh^{-\frac{3p}{2}}
           +
           K\big( \E[|X_0|^{\bar{p}}] \big)^{\frac{p}{\bar{p}}}
      \big)
      \big(
           1
           +
           \big(\E[|X_{0}|^{\bar{p}}]\big)^{\frac{1}{q'}}
      \big)
      \\\leq&
      K
      \big(
           1
           +
           \big(\E[|X_{0}|^{\bar{p}}]\big)^{\frac{p}{\bar{p}}+\frac{1}{q'}}
      \big).
\end{split}
\end{equation*}
This together with \eqref{eq:tame.parta} implies
\begin{equation*}
\begin{split}
      &\E\big[|Y_n|^{p}\big]
      =
      \E\big[ \1_{\Omega_{R(h),n}}|Y_n|^{p} \big]
      +
      \E\big[ \1_{\Omega_{R(h),n}^c}|Y_n|^{p} \big]
      \\\leq&
      \big(\E[ \1_{\Omega_{R(h),n}}|Y_n|^{\bar{p}} ]\big)^{\frac{p}{\bar{p}}}
      +
      \E\big[ \1_{\Omega_{R(h),n}^c}|Y_n|^{p} \big]
      \leq
      K
      \big(
           1
           +
           \big(\E[|X_{0}|^{\bar{p}}]\big)^{\frac{p}{\bar{p}}+\frac{1}{q'}}
      \big)
\end{split}
\end{equation*}
for all $p \in \big[2,\tfrac{\bar{p}-H(q)}{1+\frac{3}{2}H(q)}\big]$,
which immediately yields \eqref{eq:tame.monent.bound}.
\end{proof}

\subsection{Convergence rates of the tamed Euler method}
\label{subsection:conv-rate-tamed-Euler}

Here we will detect the local convergence rates $p_1$ and $p_2$ from \eqref{eq:ed}--\eqref{eq:msd} and thus derive the global convergence rates of the tamed Euler method \eqref{eq:tame.method} via
Theorem \ref{theorem:convergence}. 

\subsubsection{Convergence rates under polynomial growth condition}

\begin{theorem}\label{Theorem:tame.order.I}
      Suppose Assumptions \ref{as:assumption1}, \ref{as:assumption2}, \ref{ass:drift-polynomial} hold  and Let $\{X(t)\}_{0 \leq t \leq T}$ and $\{Y_n\}_{0\leq n \leq N}$ be given by \eqref{eq:SIDE} and \eqref{eq:tame.method}, respectively. Let $H(q)=\max\left\{ 1+q,\tfrac{3}{2}q \right\}$ and let $\bar{p} \geq 2+4H(q)$ be a sufficiently large even number. Also let $\kappa > 0$ satisfy $q\leq\kappa\bar{p}$, $(2+q)/(1-\kappa)\leq\bar{p}$ with $\bar{p}\geq2+2q$.
      Then
      there exist $\gamma > 0$ and $K > 0$ independent of $h$ such that
      \begin{equation*}
            \sup_{ 0 \leq n \leq N }
            \big( \E\big[ |X(t_{n}) - Y_{n}|^2 \big] \big)^{\frac{1}{2}}
            \leq
            K\big( 1+(\E[|X_{0}|^{\bar{p}}])^{\gamma} \big) h^{\frac{1}{2}-\kappa}.
      \end{equation*}
\end{theorem}

\begin{proof}
Consider the one-step approximation of \eqref{eq:tame.method}
\begin{equation}\label{TMethodVersion.sub}
\begin{split}
      Y_{t,x}(t+h)
      =
      x&+\frac{f(x)h}{1+|f(x)|h}
      +
      \frac{g(x)(W(t+h)-W(t))}{1+|g(x)|h}
      \\&+
      \int_{t}^{t+h}\int_Z
          \frac{\sigma(x,z)}{1+|\sigma(x,z)|h}
      \,\bar{N}(\mbox{d}s,\mbox{d}z)
\end{split}
\end{equation}
and the one-step approximation of the Euler-Maruyama method
\begin{equation}\label{eq:em.sub}
      Y^E_{t,x}(t+h)
      =
      x+f(x)h+g(x)(W(t+h)-W(t))
      +
      \int_{t}^{t+h}\int_Z\sigma(x,z)\,\bar{N}(\mbox{d}s,\mbox{d}z).
\end{equation}
To discuss $\E[X_{t,x}(t+h)-Y_{t,x}(t+h)]$, we decompose it as follows
\begin{equation}\label{eq:expectation.sub}
\begin{split}
      \big|\E[X_{t,x}(t+h)-Y_{t,x}(t+h)]\big|
      \leq&
      \big|\E[X_{t,x}(t+h)-Y^E_{t,x}(t+h)]\big|
      \\&+
      \big|\E[Y^E_{t,x}(t+h)-Y_{t,x}(t+h)]\big|.
\end{split}
\end{equation}
It follows from \eqref{eq:SIDE}, \eqref{eq:em.sub} and the martingale property that
\begin{equation}\label{eq:tame.emexpect.sub}
\begin{split}
      \big| \E[X_{t,x}(t+h)-Y^E_{t,x}(t+h)] \big|
      \leq
      \int_t^{t+h}\E\big[|f(X_{t,x}(s^{-}))-f(x)|\big]\,\mbox{d}s.
\end{split}
\end{equation}
We then apply \eqref{eq:fpg}, \eqref{eq:anal.solu.bounded} with $\kappa \geq \frac{q}{\bar{p}} \geq \frac{q}{2\bar{p}}$ and the H\"{o}lder inequality to derive
\begin{equation}\label{eq:tame.a.sub}
\begin{split}
      \E\big[|f(X_{t,x}&(s^{-}))-f(x)|\big]
      \leq
      K\E\big[
         \big(1+|X_{t,x}(s^{-})|^{\frac{q}{2}}+|x|^{\frac{q}{2}}\big)
         |X_{t,x}(s^{-})-x|
         \big]
      \\\leq&
        K\big(
         \E\big[
           \big(
               1
               +
               |X_{t,x}(s^{-})|^{\frac{q}{2}}
               +
               |x|^{\frac{q}{2}}
           \big)^{\frac{1}{\kappa}}
           \big]
         \big)^{\kappa}
         \big(
         \E\big[
         |X_{t,x}(s^{-})-x|^{\frac{1}{1-\kappa}}
         \big]
         \big)^{1-\kappa}
      \\\leq&
      K\big(1+|x|^{\frac{q}{2}}\big)
       \big(
         \E\big[
           |X_{t,x}(s^{-})-x|^{\frac{1}{1-\kappa}}
           \big]
       \big)^{1-\kappa}.
\end{split}
\end{equation}
We use \eqref{eq:SIDE}, \eqref{eq:elementaryinequality}, the H\"{o}lder inequality and \eqref{eq:bdg} with $\frac{1}{1-\kappa} \downarrow 1$ to get
\begin{equation}\label{eq:tame.b.sub}
\begin{split}
      \E\big[|X_{t,x}(s^{-})&-x|^{\frac{1}{1-\kappa}}\big]
      \leq
      K\E\Big[
         \Big|\int_{t}^{s}f(X_{t,x}(r^{-}))\,\mbox{d}r\Big|^{\frac{1}{1-\kappa}}
         \Big]
      \\&+
      K\Big(\E\Big[
         \Big|\int_{t}^{s}g(X_{t,x}(r^{-}))\,\mbox{d}W(r)\Big|^{\frac{2}{1-\kappa}}
         \Big]\Big)^{\frac{1}{2}}
      \\&+
      K\E\Big[
         \Big|\int_{t}^{s}\int_{Z}\sigma(X_{t,x}(r^{-}),z)\,\bar{N}(\mbox{d}r,\mbox{d}z)\Big|^{\frac{1}{1-\kappa}}
         \Big]
      \\\leq&
            K(s-t)^{\frac{1}{1-\kappa}-1}
            \int_{t}^{s}
                \E\big[|f(X_{t,x}(r^{-}))|^{\frac{1}{1-\kappa}}\big]
            \,\mbox{d}r
            \\&+
            K\Big((s-t)^{\frac{1}{1-\kappa}-1}
             \int_{t}^{s}
                 \E\big[ |g(X_{t,x}(r^{-}))|^{\frac{2}{1-\kappa}} \big]
             \,\mbox{d}r
            \Big)^{\frac{1}{2}}
            \\&+
            K\E\Big[
            \int_{t}^{s}\int_{Z}
                \big|\sigma(X_{t,x}(r^{-}),z)\big|^{2}
             \,\nu(\mbox{d}z)\mbox{d}r
         \Big]^{\frac{1}{2(1-\kappa)}}.
\end{split}
\end{equation}
Then H\"older's inequality, \eqref{eq:fsg}--\eqref{eq:hsg} and \eqref{eq:anal.solu.bounded} with $\bar{p} \geq \frac{2+q}{2(1-\kappa)} \geq \frac{2+q}{2(1-\kappa)}$ and $\bar{p} \geq 2+\frac{q}{2}$ promise
      \begin{equation}\label{eq:tame.bb.sub}
      \begin{split}
            \E\big[| X_{t,x}(s) &- x |^{\frac{1}{1-\kappa}}\big]
            \leq
            K(s-t)^{\frac{1}{1-\kappa}}
            \big(1+|x|^{\frac{2+q}{2(1-\kappa)}}\big)
            \\&+
            K(s-t)^{\frac{1}{2(1-\kappa)}}
            \big(1+|x|^{\frac{4+q}{4(1-\kappa)}}\big)
      \leq
             Kh^{\frac{1}{2(1-\kappa)}}\big(1+|x|^{\frac{2+q}{2(1-\kappa)}}\big).
      \end{split}
      \end{equation}
A combination of \eqref{eq:tame.bb.sub}, \eqref{eq:tame.emexpect.sub} and \eqref{eq:tame.a.sub} gives
\begin{equation}\label{eq:tame.emexpect.a.sub}
      \big| \E[X_{t,x}(t+h)-Y^E_{t,x}(t+h)] \big|
      \leq
      Kh^{\frac{3}{2}}\big(1+|x|^{1+q}\big).
\end{equation}
By \eqref{TMethodVersion.sub}--\eqref{eq:em.sub} and \eqref{eq:fsg}, it is easy to see that
\begin{equation}\label{eq:LF1.sub}
\begin{split}
      \big|\E[ Y^E_{t,x}(t+h)-Y_{t,x}(t+h) ]\big|
      =
      \frac{|f(x)|^2h^2}{1+|f(x)|h}
      \leq
      Kh^{2}\big( 1+|x|^{2+q} \big).
\end{split}
\end{equation}
Substituting \eqref{eq:tame.emexpect.sub} and \eqref{eq:LF1.sub} into \eqref{eq:expectation.sub} shows \eqref{eq:ed} is satisfied with $p_{1} = \frac{3}{2}$.
Next we examine the one-step error in mean-square sense.
By \eqref{eq:elementaryinequality}, we have
\begin{equation}\label{Deviation.tamed.sub}
\begin{split}
      \E\big[ |X_{t,x}(t+h)-Y_{t,x}(t+h)|^{2} \big]
      \leq&
      2\E\big[ |X_{t,x}(t+h)-Y^E_{t,x}(t+h)|^{2} \big]
      \\&+
      2\E\big[ |Y^E_{t,x}(t+h)-Y_{t,x}(t+h)|^{2} \big].
\end{split}
\end{equation}
It follows from \eqref{eq:SIDE}, \eqref{eq:em.sub}, the H\"{o}lder inequality and isometry formulae that
      \begin{equation}\label{EulerMean-SquareE.sub.a}
      \begin{split}
            \E\big[|X_{t,x}(t+h)&-Y^E_{t,x}(t+h)|^{2}\big]
            \leq
            3h\int_t^{t+h}
                  \E\big[|f(X_{t,x}(s^{-}))-f(x)|^{2}\big]
              \,\mbox{d}s
            \\&+
            3\int_t^{t+h}
                 \E[|g(X_{t,x}(s^{-}))-g(x)|^{2}]
             \,\mbox{d}s
            \\&+
            3\int_t^{t+h}
                 \E\Big[
                     \int_Z
                         |\sigma(X_{t,x}(s^{-}),z)-\sigma(x,z)|^{2}
                     \,\nu(\mbox{d}z)
                     \Big]
             \mbox{d}s.
      \end{split}
      \end{equation}
      Applying techniques used in \eqref{eq:tame.a.sub} yields
      \begin{equation}\label{eq:subsubsubsub}
      \begin{split}
            \E\big[|f(X_{t,x}(s))-f(x)|^{2}\big]
      \leq&
            K
            \big( 1+|x|^{q} \big)
            \big(
            \E\big[| X_{t,x}(s) - x |^{\frac{2}{1-\kappa}}\big]
            \big)^{1-\kappa}.
      \end{split}
      \end{equation}
      Similarly to \eqref{eq:tame.b.sub}--\eqref{eq:tame.bb.sub} and noting $\frac{1}{1-\kappa} \downarrow 1$, we can derive that
      \begin{equation}\label{eq:f.mean.square.estimation.one.sub}
      \begin{split}
            \E\big[| &X_{t,x}(s) - x |^{\frac{2}{1-\kappa}}\big]
            \leq
            K(s-t)^{\frac{2}{1-\kappa}-1}
            \int_{t}^{s}
                \E\big[|f(X_{t,x}(r^{-}))|^{\frac{2}{1-\kappa}}\big]
            \,\mbox{d}r
            \\&+
            K\big(
                  1
                  +
                  (\nu(Z)(s-t))^{\frac{1}{1-\kappa}-1}
            \big)
            \int_{t}^{s}
                 \E\Big[
                   \int_{Z}
                       |\sigma(X_{t,x}(r^{-}),z)|^{\frac{2}{1-\kappa}}
                   \,\nu(\mbox{d}z)
                 \Big]
            \mbox{d}r
            \\&+
            K(s-t)^{\frac{1}{1-\kappa}-1}
            \int_{t}^{s}
                \E\big[|g(X_{t,x}(r^{-}))|^{\frac{2}{1-\kappa}}\big]
            \,\mbox{d}r
      \leq
             Kh\big(1+|x|^{\frac{2+q}{1-\kappa}}\big).
      \end{split}
      \end{equation}
      Inserting \eqref{eq:f.mean.square.estimation.one.sub} into \eqref{eq:subsubsubsub} gives
      \begin{equation}\label{eq:fsubsubsubsub}
            \E\big[|f(X_{t,x}(s))-f(x)|^{2}\big]
      \leq
            Kh^{1-\kappa}\big(1+|x|^{2+2q}\big).
      \end{equation}
      Likewise, one can prove
      \begin{align}
            \label{eq:gsubsubsubsub}
            \E\big[|g(X_{t,x}(s))-g(x)|^{2}\big]
            \leq&
            Kh^{1-\kappa}\big(1+|x|^{2+2q}\big),
            \\
            \label{eq:sigmasubsubsubsub}
            \E\Big[
              \int_{Z}
                  |\sigma(X_{t,x}(s^{-}),z)-\sigma(x,z)|^{2}
              \,\nu(\mbox{d}z)
              \Big]
            \leq&
            Kh^{1-\kappa}\big(1+|x|^{2+2q}\big).
      \end{align}
      Then \eqref{eq:fsubsubsubsub}--\eqref{eq:sigmasubsubsubsub} and \eqref{EulerMean-SquareE.sub.a} enable us to obtain
      \begin{equation}\label{EulerMean-SquareE.sub.b}
            \E\big[|X_{t,x}(t+h)-Y^E_{t,x}(t+h)|^{2}\big]
            \leq
            Kh^{2-\kappa}\big(1+|x|^{2+2q}\big).
      \end{equation}
Moreover, by \eqref{TMethodVersion.sub}--\eqref{eq:em.sub} and  \eqref{eq:fsg}--\eqref{eq:hsg} we derive
\begin{equation}\label{LF2tamed.sub}
\begin{split}
      \E\big[| Y^E_{t,x}(t+h)&-Y_{t,x}(t+h) |^{2}\big]
      \leq
      3\Big|f(x)h\frac{|f(x)|h}{1+|f(x)|h}\Big|^{2}
      \\&+
      3\E\Big[
         \Big|g(x)(W(t+h)-W(t))\frac{|g(x)|h}{1+|g(x)|h}\Big|^2
         \Big]
      \\&+
      3\E\Big[
         \Big|
           \int_t^{t+h}\int_Z
               \sigma(x,z)\frac{|\sigma(x,z)|h}{1+|\sigma(x,z)|h}
           \,\bar{N}(\mbox{d}s,\mbox{d}z)\Big|^2\Big]
      \\\leq&
      3h^{4}|f(x)|^{4}+3h^{2}|g(x)|^{4}\E\big[|W(t+h)-W(t)|^{2}\big]
      \\&+
      3h^{2}
      \int_{t}^{t+h}\int_{Z}
          |\sigma(x,z)|^{4}
      \,\nu(\mbox{d}z)\,\mbox{d}s
      \leq
      Kh^{3} \big(1+|x|^{4+2q}\big).
\end{split}
\end{equation}
Plugging \eqref{EulerMean-SquareE.sub.b}--\eqref{LF2tamed.sub} into \eqref{Deviation.tamed.sub} implies \eqref{eq:msd} is satisfied with $p_{2}=1-\kappa$. Finally, Theorem \ref{theorem:convergence} gives the desired order and finishes the proof. \end{proof}

At this moment, we would like to point out that the mean-square convergence order of the tamed Euler method \eqref{eq:tame.method}, arbitrarily close to $\frac{1}{2}$, coincides with that in \cite[Theorem 3.5]{Dareiotis16} and
\cite[Theorem 2]{kumar2017explicit}, covering a wider class of L\'{e}vy noise.
Different from \cite{Dareiotis16, kumar2017explicit}, we allow jump coefficients to grow super-linearly but require finite L\'{e}vy measure.
%
When $\sigma \equiv 0$, i.e., the jump-diffusion SDEs \eqref{eq:JSDE} reduce to the continuous SDEs and
the corresponding numerical results of such equation in \cite{andersson2016mean,Hutzenthaler15,Hutzenthaler12,Mao13a,
Sabanis16,sabanis2013note,Tretyakov13} can be recovered.

\begin{corollary}\label{Theorem:nonjump}
      Suppose Assumptions \ref{as:assumption1}, \ref{as:assumption2} and \ref{ass:drift-polynomial} with $\sigma \equiv 0$ hold and let $\{X(t)\}_{0 \leq t \leq T}$ and $\{Y_n\}_{0\leq n \leq N}$ be given by \eqref{eq:SIDE} and \eqref{eq:tame.method}, respectively.
      Let $H(q)=\max\left\{ 1+q,\tfrac{3}{2}q \right\}$ and let $\bar{p} \geq \max\{2+4H(q),2+2q\}$ be a sufficiently large even number.
      Then there exist $\gamma > 0$ and $K > 0$ independent of $h$ such that
      \begin{equation*}
            \sup_{ 0 \leq n \leq N }
            \big( \E\big[ |X(t_{n}) - Y_{n}|^2 \big] \big)^{\frac{1}{2}}
            \leq
            K\big( 1+\big(\E[|X_{0}|^{\bar{p}}]\big)^{\gamma} \big) h^{\frac{1}{2}}.
      \end{equation*}
\end{corollary}

\begin{proof}
Since Theorem \ref{Theorem:tame.order.I} shows $p_1=\frac{3}{2}$, it suffices to prove $p_2=1$. To this end, we need re-eveulate $\E\big[|X_{t,x}(t+h)-Y^E_{t,x}(t+h)|^{2}\big]$. By \eqref{eq:f.mean.square.estimation.one.sub}, \eqref{eq:fsg}--\eqref{eq:gsg} and \eqref{eq:anal.solu.bounded}, we have
      \begin{equation*}
      \begin{split}
            \E\big[|X_{t,x}(s) &- x |^{\frac{2}{1-\kappa}}\big]
            \leq
            Kh^{\frac{2}{1-\kappa}-1}
            \int_{t}^{s}
                \E\big[|f(X_{t,x}(r^{-}))|^{\frac{2}{1-\kappa}}\big]
            \,\mbox{d}r
            \\&+
            Kh^{\frac{1}{1-\kappa}-1}
            \int_{t}^{s}
                \E\big[|g(X_{t,x}(r^{-}))|^{\frac{2}{1-\kappa}}\big]
            \,\mbox{d}r
      \leq
             Kh^{\frac{1}{1-\kappa}}
             \big(1+|x|^{\frac{2+q}{1-\kappa}}\big),
      \end{split}
      \end{equation*}
      where $\kappa$ is the same as that in \eqref{eq:f.mean.square.estimation.one.sub}. This and \eqref{eq:subsubsubsub} imply
      \begin{equation}\label{eq:subsubsubsub.nonjumpf}
            \E\big[|f(X_{t,x}(s))-f(x)|^{2}\big]
      \leq
            Kh\big( 1+|x|^{2+2q} \big).
      \end{equation}
      Similarly, we get
      \begin{equation}\label{eq:subsubsubsub.nonjumpg}
            \E\big[|g(X_{t,x}(s))-g(x)|^{2}\big]
            \leq
            Kh\big(1+|x|^{2+2q}\big).
      \end{equation}
      Applying \eqref{eq:subsubsubsub.nonjumpf}--\eqref{eq:subsubsubsub.nonjumpg} leads to
      \begin{equation*}\label{EulerMean.sub.a.nonjump}
      \begin{split}
            \E\big[|X_{t,x}(t+h)&-Y^E_{t,x}(t+h)|^{2}\big]
            \leq
            2h\int_t^{t+h}
                  \E\big[|f(X_{t,x}(s^{-}))-f(x)|^{2}\big]
              \,\mbox{d}s
            \\&+
            2\int_t^{t+h}
                 \E\big[|g(X_{t,x}(s^{-}))-g(x)|^{2}\big]
             \,\mbox{d}s
      \leq
            Kh^{2}\big(1+|x|^{2+2q}\big),
      \end{split}
      \end{equation*}
      which together with \eqref{LF2tamed.sub} yields $p_{2}=1$ and thus ends the proof by Theorem \ref{theorem:convergence}.
\end{proof}

\subsubsection{Higher convergence rate in the additive noise case}


We will further investigate convergence rate of \eqref{eq:tame.method} for jump-diffusion SDEs with additive noise under the following assumption.

\begin{assumption}\label{as:assumption3}
      Assume that for all $i,j,k=1,\ldots,d$, the derivatives of the coefficients $f(x) = (f_i)_{d \times 1}$ in \eqref{eq:JSDE}, i.e.,
      $
         \frac{\partial{f_{i}}}{\partial{x_{j}}},
         \frac{\partial^{2}{f_{i}}}{\partial{x_{j}}\partial{x_{k}}}$,
      are continuous and satisfy the polynomial growth condition in the form of \eqref{eq:fpg}, i.e., there exist $K,q \geq 0$ such that for all $x,y\in\R^d$, $i,j,k=1,\ldots,d$,
      \begin{equation*}
            |a(x)-a(y)|^{2}
            \leq
            K(1+|x|^{q}+|y|^{q})|x-y|^{2},
            \quad
            a
            :=
            \frac{\partial{f_{i}}}{\partial{x_{j}}},
            \frac{\partial^{2}{f_{i}}}{\partial{x_{j}}\partial{x_{k}}}.
      \end{equation*}
\end{assumption}

The following theorem, in some sense, can be regarded as an extension of existing
known results for the additive noise case in our setting.

\begin{theorem}\label{Theorem:tame.additive.order}
      Suppose Assumptions \ref{as:assumption1}, \ref{as:assumption2},  \ref{as:assumption3} with $g(x)=g \in \R^{d{\times}m}, \sigma(x,z)=\sigma(z) \in \R^d$ for all $x\in\R^d, z \in Z$ hold and let $\{X(t)\}_{0 \leq t \leq T}$ and $\{Y_n\}_{0\leq n \leq N}$ be given by \eqref{eq:SIDE} and \eqref{eq:tame.method}, respectively. Let $H(q)=\max\left\{ 1+q,\tfrac{3}{2}q \right\}$ and let $\bar{p} \geq 2+4H(q)$ be a sufficiently large even number.
      Then there exist $\gamma > 0$ and $K > 0$ independent of $h$ such that
      \begin{equation*}
            \sup_{ 0 \leq n \leq N }
            \big( \E\big[ |X(t_{n}) - Y_{n}|^2 \big] \big)^{\frac{1}{2}}
            \leq
            K\big( 1+(\E[|X_{0}|^{\bar{p}}])^{\gamma} \big) h.
      \end{equation*}
\end{theorem}

\begin{proof}
Using Lemma \ref{lemma:itoformula} shows that for all $i=1,\ldots,d$,
\begin{equation*}
\begin{split}
      f_{i}(X_{t,x}(s))&-f_{i}(x)
      =
      \int_{t}^{s}
          \sum_{j=1}^{d}
          \frac{\partial{f_{i}}(X_{t,x}(r^{-}))}
                {\partial{x_{j}}}
          f_{j}(X_{t,x}(r^{-}))
      \,\mbox{d}r
      \\&+
      \frac{1}{2}
      \int_{t}^{s}
          \sum_{j=1}^{d}\sum_{k=1}^{d}\sum_{l=1}^{m}
          \frac{\partial^{2}{f_{i}}(X_{t,x}(r^{-}))}
                {\partial{x_{j}}\partial{x_{k}}}
          g_{j,l}g_{k,l}
      \,\mbox{d}r
      \\&+
      \int_{t}^{s}\int_{Z}
          f_{i}
          \big(X_{t,x}(r^{-})+\sigma(z)\big)-f_{i}(X_{t,x}(r^{-}))
          \\&-
          \sum_{j=1}^{d}
          \frac{\partial{f_{i}}(X_{t,x}(r^{-}))}
                {\partial{x_{j}}}
               \sigma_{j}(z)
      \,\nu(\mbox{d}z)\mbox{d}r
      \\&+
      \sum_{k=1}^{m}
      \int_{t}^{s}
          \sum_{j=1}^{d}
          \frac{\partial{f_{i}}(X_{t,x}(r^{-}))}
                {\partial{x_{j}}}
          g_{j,k}
      \,\mbox{d}W_{k}(r)
      \\&+
      \int_{t}^{s}\int_{Z}
          f_{i}
          \big(X_{t,x}(r^{-})+\sigma(z)\big)-f_{i}(X_{t,x}(r^{-}))
      \,\bar{N}(\mbox{d}r,\mbox{d}z)
      \\:=& B_{1}+B_{2}+B_{3}+B_{4}+B_{5},
        ~\P\text{-a.s.},
\end{split}
\end{equation*}
which together with the martingale property yields
\begin{equation*}
      \big| \E[f(X_{t,x}(s))-f(x)] \big|^{2}
      =
      \sum\limits_{i=1}^{d}
      \big| \E[f_{i}(X_{t,x}(s))-f_{i}(x)] \big|^{2}
      =
      \sum\limits_{i=1}^{d}
      \big| \E[B_{1} + B_{2} + B_{3}] \big|^{2}.
\end{equation*}
By \eqref{eq:elementaryinequality}, Assumption \ref{as:assumption3} and \eqref{eq:anal.solu.bounded}, we have
\begin{equation*}
\begin{split}
      \big| \E&[f(X_{t,x}(s))-f(x)] \big|^{2}
      \leq
      3\sum\limits_{i=1}^{d}
      \Big(
      \int_{t}^{s}
          \E\Big[\Big|
          \sum_{j=1}^{d}
          \frac{\partial{f_{i}}(X_{t,x}(r^{-}))}
                {\partial{x_{j}}}
          f_{j}(X_{t,x}(r^{-}))
          \Big|\Big]
      \,\mbox{d}r
      \Big)^{2}
      \\&+
      \Big(\int_{t}^{s}
          \E\Big[\Big|\int_{Z}
          f_{i}
          \big(X_{t,x}(r^{-})+\sigma(z)\big)-f_{i}(X_{t,x}(r^{-}))
      \\&-
      \sum_{j=1}^{d}
          \frac{\partial{f_{i}}(X_{t,x}(r^{-}))}
                {\partial{x_{j}}}
                \sigma_{j}(z)
      \,\nu(\mbox{d}z)
      \Big|\Big]\mbox{d}r
      \Big)^{2}
      \\&+
      \Big(\frac{1}{2}
      \int_{t}^{s}
          \E\Big[\Big|
          \sum_{j=1}^{d}\sum_{k=1}^{d}\sum_{l=1}^{m}
          \frac{\partial^{2}{f_{i}}(X_{t,x}(r^{-}))}
                {\partial{x_{j}}\partial{x_{k}}}
          g_{j,l}g_{k,l}
          \Big|\Big]
      \,\mbox{d}r
      \Big)^{2}
      \leq
      K(s-t)^{2}\big(1+|x|^{2+q}\big)^{2}.
\end{split}
\end{equation*}
This and \eqref{eq:tame.emexpect.sub} imply
\begin{equation}\label{eq:tame.emexpect.additive}
\begin{split}
      \big| \E\big[ X_{t,x}(t+h)-Y^{E}_{t,x}(t+h) \big] \big|
      \leq&
      \int_{t}^{t+h}
          \big| \E[ f(X_{t,x}(s^{-}))-f(x) ] \big|
      \,\mbox{d}s
      \\\leq&
      Kh^{2}(1+|x|^{2+q}).
\end{split}
\end{equation}
Similarly to \eqref{eq:LF1.sub}, we can obtain
\begin{equation*}\label{eq:LF1A}
\begin{split}
      \big| \E\big[ Y^{E}_{t,x}(t+h)-Y_{t,x}(t+h) \big] \big|
      \leq
      Kh^{2}\big(1+|x|^{2+q}\big).
\end{split}
\end{equation*}
The triangle inequality suggests that \eqref{eq:ed} is satisfied with $p_1=2$.
Thanks to \eqref{eq:elementaryinequality},
\begin{equation}\label{eq:1}
\begin{split}
      &\E\big[ |f_{i}(X_{t,x}(s))-f_{i}(x)|^{2} \big]
      \\\leq&
      5\E\big[|B_{1}|^{2}\big]
      +
      5\E\big[|B_{2}|^{2}\big]
      +
      5\E\big[|B_{3}|^{2}\big]
      +
      5\E\big[|B_{4}|^{2}\big]
      +
      5\E\big[|B_{5}|^{2}\big].
\end{split}
\end{equation}
Since the first three terms on the right hand side of \eqref{eq:1} can be estimated in the same manner, here we just, for example, give the estimate of $5\E\big[|B_{2}|^{2}\big]$ via the H\"older inequality, Assumption \ref{as:assumption3} and \eqref{eq:anal.solu.bounded} as follows
\begin{equation*}\label{eq:2}
\begin{split}
      5\E\big[|B_{2}|^{2}\big]
      \leq&
      \frac{5}{4}(s-t)
      \int_{t}^{s}
          \E\bigg[
            \Big|
            \sum_{j=1}^{d}\sum_{k=1}^{d}\sum_{l=1}^{m}
            \frac{\partial^{2}{f_{i}}(X_{t,x}(r^{-}))}
                 {\partial{x_{j}}\partial{x_{k}}}
            g_{j,l}g_{k,l}
            \Big|^{2}
            \bigg]
      \,\mbox{d}r
      \\=&
      \frac{5}{4}(s-t)
      \sum_{j,j'=1}^{d}\sum_{k,k'=1}^{d}\sum_{l,l'=1}^{m}
      \int_{t}^{s}
          \E\bigg[
            \Big|
            \frac{\partial^{2}{f_{i}}(X_{t,x}(r^{-}))}
                 {\partial{x_{j}}\partial{x_{k}}}
            \frac{\partial^{2}{f_{i}}(X_{t,x}(r^{-}))}
                 {\partial{x_{j'}}\partial{x_{k'}}}
                 \\&g_{j,l}g_{k,l}g_{j',l'}g_{k',l'}
            \Big|
            \bigg]
      \,\mbox{d}r
      \\\leq&
      K(s-t)\int_{t}^{s}
            \big(1+\E\big[|X_{t,x}(r^{-})|^{2+q}\big]\big)
        \,\mbox{d}r
      \leq
      K(s-t)^{2}(1+|x|^{2+q}).
\end{split}
\end{equation*}
Similarly, $5\E\big[|B_{4}|^{2}\big]$ and $5\E\big[|B_{5}|^{2}\big]$ are calculated by
\begin{equation*}
\begin{split}
      5\E\big[|B_{4}|^{2}\big]
      =&
      5\sum_{k=1}^{m}
       \int_{t}^{s}
           \E\bigg[
             \Big|
             \sum_{j=1}^{d}
             \frac{\partial{f_{i}}(X_{t,x}(r^{-}))}
                  {\partial{x_{j}}} g_{j,k}
             \Big|^{2}
             \bigg]
       \,\mbox{d}r
      \\=&
      5\sum_{k=1}^{m}
       \sum_{j,j'=1}^{d}
       \int_{t}^{s}
           \E\bigg[
             \frac{\partial{f_{i}}(X_{t,x}(r^{-}))}
                  {\partial{x_{j}}}
             \frac{\partial{f_{i}}(X_{t,x}(r^{-}))}
                  {\partial{x_{j'}}} g_{j,k} g_{j',k}
             \bigg]
       \,\mbox{d}r
      \\\leq&
      K\int_{t}^{s}
            \big(
               1
               +
               \E\big[|X_{t,x}(r^{-})|^{2+q}\big]
            \big)\,\mbox{d}r
      \leq
      K(s-t)(1+|x|^{2+q})
\end{split}
\end{equation*}
and
\begin{equation*}\label{eq:3}
\begin{split}
      5\E\big[|B_{5}|^{2}\big]
      =&
      5\E\bigg[
         \int_{t}^{s}\int_{Z}
             \big|
             f_{i}\big(X_{t,x}(r^{-})+\sigma(z)\big)
             -
             f_{i}(X_{t,x}(r^{-}))
             \big|^{2}
         \,\nu(\mbox{d}z)\mbox{d}r
         \bigg]
      \\\leq&
      K\E\bigg[
         \int_{t}^{s}\int_{Z}
             \big(
                  1
                  +
                  |\sigma(z)|^{q}
                  +
                  |X_{t,x}(r^{-})|^{q}
             \big)
             |\sigma(z)|^{2}
         \,\nu(\mbox{d}z)\mbox{d}r
         \bigg]
      \\\leq&
      K\int_{t}^{s}
           \big(
              1+\E\big[|X_{t,x}(r^{-})|^{q}\big]
           \big)
       \mbox{d}r
      \leq
      K(s-t)\big( 1+|x|^{q} \big).
\end{split}
\end{equation*}
Combining the above estimates promises
\begin{equation*}\label{fEulerMean-SquareA}
      \E\big[|f(X_{t,x}(s))-f(x)|^{2}\big]
      =
      \E\Big[
        \sum_{i=1}^{d}|f_{i}(X_{t,x}(s))-f_{i}(x)|^{2}
        \Big]
      \leq
      K(s-t)\big(1+|x|^{2+q}\big),
\end{equation*}
which together with the H\"{o}lder inequality realizes that
\begin{equation}\label{EulerMean-SquareE.additive}
\begin{split}
      \E\big[ |X_{t,x}(t+h)-Y^{E}_{t,x}(t+h)|^{2} \big]
      =&
      \E\Big[
        \Big|\int_{t}^{t+h}f(X_{t,x}(s^{-}))-f(x)\,\mbox{d}s\Big|^{2}
        \Big]
      \\\leq&
        h\int_{t}^{t+h}
            \E\big[ |f(X_{t,x}(s^{-}))-f(x)|^{2} \big]
        \,\mbox{d}s
      \\\leq&
      Kh^{3}\big(1+|x|^{2+q}\big).
\end{split}
\end{equation}
Similarly to \eqref{LF2tamed.sub}, we can get
\begin{equation}\label{LF2A}
      \E\big[ |Y^{E}_{t,x}(t+h)-Y_{t,x}(t+h)|^{2} \big]
      \leq
      Kh^{3}\big(1+|x|^{4+2q}\big).
\end{equation}
Combining \eqref{EulerMean-SquareE.additive}--\eqref{LF2A} and \eqref{eq:elementaryinequality} shows that \eqref{eq:msd} is satisfied with $p_2=\frac{3}{2}$, which completes the proof by Theorem \ref{theorem:convergence}.  \end{proof}

\section{Application of the fundamental convergence theorem:
         convergence rates of the sine Euler method}\label{sec:sine}

Motivated by the explicit schemes introduced in \cite{ZM17,Zhang14}, we propose the sine Euler method for \eqref{eq:JSDE}, given by $Y_0=X_0$ and
\begin{equation}\label{eq:sin}
\begin{split}
      Y_{n+1}
      =
      Y_n
      &+
      \sin(f(Y_n)h)
      +
      \frac{\sin(g(Y_n)h)}{h} \Delta{W_n}
      \\&+
      \int_{t_n}^{t_{n+1}}\hspace{-0.25em}\int_Z
          \frac{\sin(\sigma(Y_n,z)h)}{h}
      \,\bar{N}(\mbox{d}s,\mbox{d}z),
      \quad n=0,1,\ldots,N-1,
\end{split}
\end{equation}
where $\sin(x):=(\sin(x_{i}))_{d \times 1}, \forall x \in \R^d$ and $\sin(y):=(\sin(y_{ij}))_{d \times m}, \forall y \in \R^{d \times m}$.
Scheme \eqref{eq:sin} is different from schemes in \cite{ZM17,Zhang14} even if the jump term vanishes.

\subsection{Bounded $p$-th moments of the sine Euler method}
\label{subsec:bonded-mom}

\begin{lemma}\label{lemma:sinebound}
      Suppose Assumptions \ref{as:assumption1}, \ref{as:assumption2}, \ref{ass:drift-polynomial} hold and let $\{Y_n\}_{0 \leq n \leq N}$ be given by \eqref{eq:sin}. Let $H(q)=\max\left\{ 1+q,\tfrac{3}{2}q \right\}$ and let $\bar{p} \geq 2+4H(q)$ be a sufficiently large even number. Then there exist $\beta > 0$ and $K > 0$ independent of $h$ such that
      \begin{equation*}\label{eq:sinebound}
            \sup\limits_{0 \leq n \leq N }
            \E\big[|Y_n|^p\big]
            \leq
            K\big( 1+\big(\E[|X_0|^{\bar{p}}]\big)^{\beta} \big),
            \quad
            \forall p \in \big[2,\tfrac{\bar{p}-H(q)}{1+\frac{3}{2}H(q)}\big].
      \end{equation*}
\end{lemma}

\begin{proof}
Since $|\sin z|\leq 1, \forall z \in \R$, we have $|\sin (x)| \leq \sqrt{d}, \forall x \in \R^{d}$ and $|\sin (y)| \leq \sqrt{md}, \forall y \in \R^{d \times m}$, which together with \eqref{eq:sin} gives
\begin{equation*}\label{eq:sineb}
\begin{split}
      &|Y_{n+1}|
      \leq
      |Y_n|
      +
      \sqrt{d}
      +
      \sqrt{md}h^{-1}|\Delta{W_n}|
      +
      \sqrt{d}\nu(Z)
      +
      \int_{t_n}^{t_{n+1}}\int_{Z}
          \sqrt{d}h^{-1}
      \,N(\mbox{d}s,\mbox{d}z)
      \\&\leq
      |X_0|
      +
      \sqrt{d}(1+\nu(Z))(n+1)
      +
      \sqrt{md}h^{-1}
      \sum\limits_{i=0}^{n}|\Delta{W_i}|
      +\sqrt{d}h^{-1}
      \int_{0}^{t_{n+1}}\hspace{-0.5em}\int_Z\hspace{-0.25em}
      \,N(\mbox{d}s,\mbox{d}z).
\end{split}
\end{equation*}
Let $R>0$ be sufficiently large and introduce a sequence of decreasing sets
\begin{equation*}\label{Eq:Omega}
      \Omega_{R,n}
      :=
      \{
        \omega \in \Omega
        :
        \sup_{0 \leq i \leq n}
        |Y_i(\omega)|
        \leq
        R
      \},
      \quad
      \forall n = 0,1,\ldots,N-1, N \in \N.
\end{equation*}
We define a continuous-time approximation $\{\bar{Y}(t)\}_{0\leq{t}\leq{T}}$ of $\{Y_n\}_{0\leq{n}\leq{N}}$ by
\begin{equation*}\label{eq:sine.con.c}
\begin{split}
      \bar{Y}(t)
      =
      Y_n
      &+
      \int_{t_n}^{t} \frac{\sin(f(Y_n)h)}{h} \,\mbox{d}s
      +
      \int_{t_n}^{t} \frac{\sin(g(Y_n)h)}{h} \,\mbox{d}W(s)
      \\&+
      \int_{t_n}^{t}\int_Z
          \frac{\sin(\sigma(Y_n,z)h)}{h}
      \,\bar{N}(\mbox{d}s,\mbox{d}z),
      ~~\P\text{-a.s.}
\end{split}
\end{equation*}
for all $t \in [t_n,t_{n+1}], n = 0,1,\ldots,N-1$. Similarly to Lemma \ref{lemma:tame.moment.bound.a}, we have
\begin{align*}\label{eq:sine.usmartingale}
\begin{split}
      \E\big[ &\1_{\Omega_{R, n+1}}|\bar{Y}(t)|^{\bar{p}} \big]
      \leq
      \E\big[\1_{\Omega_{R,n}}|Y_n|^{\bar{p}}\big]+\bar{p}\int_{t_n}^t\E\big[\1_{\Omega_{R,n}}\big|Y_n\big|^{\bar{p}-2}\big<Y_n,f(Y_n)\big>\big]\,\mbox{d}s
      \\&+\frac{\bar{p}(\bar{p}-1)}{2}\int_{t_n}^t\E\big[\1_{\Omega_{R,n}}\big|Y_n\big|^{\bar{p}-2}\big|g(Y_n)\big|^2\big]\,\mbox{d}s
         +K\int_{t_n}^t\E\big[\1_{\Omega_{R,n}}\big|\bar{Y}(s^-)\big|^{\bar{p}}\big]\,\mbox{d}s
      \\&+\big(1+(\bar{p}-2)\varepsilon\big)\int_{t_n}^t\E\Big[\1_{\Omega_{R,n}}\int_Z\big|\sigma(Y_n,z)\big|^{\bar{p}}\,\nu(\mbox{d}z)\Big]\mbox{d}s
         + J_{1} + J_{2} + J_{3} + J_{4},
\end{split}
\end{align*}
where
\begin{align*}
        J_{1}:=&\bar{p}\int_{t_n}^t\E\Big[\1_{\Omega_{R,n}}\big|\bar{Y}(s^-)\big|^{\bar{p}-2}\Big<\bar{Y}(s^-)-Y_n,\frac{\sin(f(Y_n)h)}{h}\Big>\Big]\,\mbox{d}s,
      \\J_{2}:=&\bar{p}\int_{t_n}^t\E\Big[\1_{\Omega_{R,n}}\big|\bar{Y}(s^-)\big|^{\bar{p}-2}\Big<Y_n,\frac{\sin(f(Y_n)h)}{h}-f(Y_n)\Big>\Big]\,\mbox{d}s,
      \\J_{3}:=&\bar{p}\int_{t_n}^t\E\big[\1_{\Omega_{R,n}}\big(|\bar{Y}(s^-)|^{\bar{p}-2}-|Y_n|^{\bar{p}-2}\big)\big<Y_n,f(Y_n)\big>\big]\,\mbox{d}s,
      \\J_{4}:=&\frac{\bar{p}(\bar{p}-1)}{2}\int_{t_n}^t\E\left[\1_{\Omega_{R,n}}\left(|\bar{Y}(s^-)|^{\bar{p}-2}-|Y_n|^{\bar{p}-2}\right)\left|g(Y_n)\right|^2\right]\,\mbox{d}s.
\end{align*}
By the coercivity condition \eqref{eq:gcc.b}, we get
\begin{equation*}\label{eq:using_gcc}
\begin{split}
      \E[ \1_{\Omega_{R,n+1}}|\bar{Y}(t)|^{\bar{p}} ]
      \leq&
      (1+Kh)
      \E[\1_{\Omega_{R,n}}|Y_{n}|^{\bar{p}}]      +
      K\int_{t_n}^{t}
           \E[\1_{\Omega_{R,n}}|\bar{Y}(s^-)|^{\bar{p}}]
       \,\mbox{d}s
      \\&+J_1+J_2+J_3+J_4.
\end{split}
\end{equation*}
Now we consider $J_1$.
Using $| \sin z | \leq |z|, \forall z \in \R$ leads to
$$
      | \sin x |
      =
      \Big(\sum\limits_{i=1}^{d} |\sin x_{i}|^{2} \Big)^{\frac{1}{2}}
      \leq
      \Big(\sum\limits_{i=1}^{d} |x_{i}|^{2} \Big)^{\frac{1}{2}}
      = |x|,
      \quad
      \forall x \in \R^{d}.
$$
This and the Schwarz inequality imply
\begin{equation*}\label{eq:I.one.a}
      J_1
      \leq
      p\int_{t_n}^{t}
           \E\big[
             \1_{\Omega_{R,n}}
             |\bar{Y}(s^-)|^{p-2}
             |\bar{Y}(s^-)-Y_{n}||f(Y_{n})|
             \big]
       \,\mbox{d}s,
\end{equation*}
further estimate of which is a copy of that of $I_1$ in \eqref{eq:I.one}.
Since $|z - \sin z| \leq |z|^{2}$ for all $z \in \R$, it holds for all $x \in \R^{d}$,
\begin{equation}\label{eq:sinineq}
           |x - \sin x |
           =
           \Big(
                \sum\limits_{i=1}^{d}
                |x_{i} - \sin x_{i}|^{2}
           \Big)^{\frac{1}{2}}
           =
           \Big(
                \sum\limits_{i=1}^{d}
                |x_{i}|^{4}
           \Big)^{\frac{1}{2}}
           \leq
           \sum\limits_{i=1}^{d}
           |x_{i}|^{2}
           =
           |x|^{2}.
\end{equation}
By the Schwarz inequality,
\begin{equation*}\label{eq:J.two}
    J_2
    \leq
    p\int_{t_n}^t
     \E\left[\1_{\Omega_{R,n}}|\bar{Y}(s^-)|^{p-2}|Y_n|\left|f(Y_n)\right|^2h\right]
     \mbox{d}s,
\end{equation*}
further estimate of which repeats that of $I_2$ in \eqref{eq:I.two}.
Additionally, $J_{3}$ and $J_{4}$ exactly coincide with $I_{3}$ and $I_{4}$, respectively.
Therefore, Lemma \ref{lemma:sinebound} is validated by repeating the proof of Lemmas \ref{lemma:tame.moment.bound.a} and \ref{lemma:tame.moment.bound}.
\end{proof}

\subsection{Convergence rates of the sine Euler method}

We analyze the convergence rates of method \eqref{eq:sin} 
 as Subsection \ref{subsection:conv-rate-tamed-Euler}
does.

\subsubsection{Convergence rates under polynomial growth condition}

\begin{theorem}\label{thm:sinerror.sub.nonadditive}
      Suppose Assumptions \ref{as:assumption1}, \ref{as:assumption2}, \ref{ass:drift-polynomial} hold and let $\{X(t)\}_{0 \leq t \leq T}$ and $\{Y_n\}_{0\leq n \leq N}$ be given by \eqref{eq:SIDE} and \eqref{eq:sin}, respectively. Let $H(q)=\max\left\{ 1+q,\tfrac{3}{2}q \right\}$ and let $\bar{p} \geq 2+4H(q)$ be a sufficiently large even number. Also let $\kappa > 0$ satisfy $q\leq\kappa\bar{p}$, $(2+q)/(1-\kappa)\leq\bar{p}$ with $\bar{p}\geq2+2q$.
      Then there exist $\gamma > 0$ and $K > 0$ independent of $h$ such that
      \begin{equation*}
            \sup_{ 0 \leq n \leq N }
            \big( \E\big[ |X(t_{n}) - Y_{n}|^2 \big] \big)^{\frac{1}{2}}
            \leq
            K\big( 1+(\E[|X_{0}|^{\bar{p}}])^{\gamma} \big)
            h^{\frac{1}{2}-\kappa}.
      \end{equation*}
\end{theorem}

\begin{proof}
We consider the one-step approximation of \eqref{eq:sin}, given by
\begin{equation}\label{eq:sin.onestep.sub}
\begin{split}
      Y_{t,x}(t+h)
      =
      x &+ \sin(f(x)h) + h^{-1}\sin(g(x)h)(W(t+h)-W(t))
      \\&+
      h^{-1}
      \int_t^{t+h}\int_{Z}
          \sin(\sigma(x,z)h)
      \,\bar{N}(\mbox{d}s,\mbox{d}z)
\end{split}
\end{equation}
and 
\eqref{eq:em.sub}.
Firstly, using \eqref{eq:sin.onestep.sub}, \eqref{eq:em.sub}, \eqref{eq:sinineq} and \eqref{eq:fsg} shows that
\begin{equation}\label{eq:sine.ed.sub}
      \big| \E[Y^E_{t,x}(t+h)-Y_{t,x}(t+h)] \big|
      =
      | f(x)h-\sin(f(x)h) |
      \leq
      Kh^{2}\big(1+|x|^{2+q}\big).
\end{equation}
We then follow arguments used in \eqref{eq:tame.emexpect.sub}--\eqref{eq:tame.emexpect.a.sub} to derive
\begin{equation}\label{eq:sine.emexpect.sub}
      \big| \E[X_{t,x}(t+h)-Y^E_{t,x}(t+h)] \big|
      \leq
      Kh^{\frac{3}{2}}\big(1+|x|^{1+q}\big).
\end{equation}
Combining \eqref{eq:sine.ed.sub} and \eqref{eq:sine.emexpect.sub},
we realize that \eqref{eq:ed} is satisfied with $p_{1} = \frac{3}{2}$.
Secondly, due to $|z-\sin z| \leq z^2$ for all $z \in \R$, we have for all $x \in \R^{d \times m}$,
\begin{equation*}\label{eq:sininequa.sub}
           |x - \sin x |
           =
           \Big(
                \sum\limits_{i=1}^{d}\sum\limits_{j=1}^{m}
                |x_{ij} - \sin x_{ij}|^{2}
           \Big)^{\frac{1}{2}}
           =
           \Big(
                \sum\limits_{i=1}^{d}\sum\limits_{j=1}^{m}
                |x_{ij}|^{4}
           \Big)^{\frac{1}{2}}
           \leq
           \sum\limits_{i=1}^{d}\sum\limits_{j=1}^{m}
           |x_{ij}|^{2}
           =
           |x|^{2}.
\end{equation*}
This together with \eqref{eq:elementaryinequality}, \eqref{eq:em.sub}, \eqref{eq:sinineq}, \eqref{eq:sin.onestep.sub} and \eqref{eq:fsg}--\eqref{eq:hsg} helps us to get
\begin{equation}\label{DeviationA.sub}
\begin{split}
      \E\big[ |Y^E_{t,x}(t+h)&-Y_{t,x}(t+h)|^{2} \big]
      \leq
      3|f(x)h-\sin(f(x)h)|^2
      \\&+
      3h^{-2}
       \E\big[
         |( g(x)h-\sin(g(x)h) )(W(t+h)-W(t)) |^2
         \big]
      \\&+
      3h^{-2}
       \E\Big[
         \Big|
           \int_t^{t+h}\int_Z
               \sigma(x,z)h-\sin(\sigma(x,z)h)
           \,\bar{N}(\mbox{d}s,\mbox{d}z)
         \Big|^{2}
         \Big]
      \\\leq&
      3h^{4}|f(x)|^{4}
      +
      3h^{3}|g(x)|^4
      +
      3h^{3}\int_{Z}|\sigma(x,z)|^{4}\,\nu(\mbox{d}z)
      \\\leq&
      Kh^3\big(1+|x|^{4+2q}\big).
\end{split}
\end{equation}
Following exactly the same lines of derivation for \eqref{EulerMean-SquareE.sub.b} guarantees
\begin{equation*}\label{eq:sine.EulerMean-SquareE.sub}
      \E\big[|X_{t,x}(t+h)-Y^E_{t,x}(t+h)|^{2}\big]
      \leq
      Kh^{2-\kappa}\big(1+|x|^{2+2q}\big),
\end{equation*}
which together with \eqref{DeviationA.sub} implies that \eqref{eq:msd} is fulfilled with $p_{2}=1-\kappa$. Thus we complete the proof by Theorem \ref{theorem:convergence}.
\end{proof}

The following result is similar to Corollary \ref{Theorem:nonjump} and its proof thus be omitted.

\begin{corollary}\label{thm:sinerror.sub.nonadditive.sigmazero}
      Suppose Assumptions \ref{as:assumption1}, \ref{as:assumption2} and \ref{ass:drift-polynomial} with $\sigma \equiv 0$ hold and let $\{X(t)\}_{0 \leq t \leq T}$ and $\{Y_n\}_{0\leq n \leq N}$ be given by \eqref{eq:SIDE} and \eqref{eq:sin}, respectively. Let $H(q)=\max\left\{ 1+q,\tfrac{3}{2}q \right\}$ and let $\bar{p} \geq \max\{2+4H(q),2+2q\}$ be a sufficiently large even number.
      Then there exist $\gamma > 0$ and $K > 0$ independent of $h$ such that
      \begin{equation*}
            \sup_{ 0 \leq n \leq N }
            \big( \E\big[ |X(t_{n}) - Y_{n}|^{2} \big] \big)^{\frac{1}{2}}
            \leq
            K\big( 1+\big(\E[|X_{0}|^{\bar{p}}]\big)^{\gamma} \big) h^{\frac{1}{2}}.
      \end{equation*}
\end{corollary}

\subsubsection{Higher convergence rate in additive noise case}
\begin{theorem}\label{thm:sinerror.additive}
      Suppose Assumptions \ref{as:assumption1}, \ref{as:assumption2} and \ref{as:assumption3} with $g(x)=g \in \R^{d{\times}m}$, $\sigma(x,z)=\sigma(z) \in \R^d$ for all $x\in\R^d, z \in Z$ hold and let $\{X(t)\}_{0 \leq t \leq T}$ and $\{Y_n\}_{0\leq n \leq N}$ be given by \eqref{eq:SIDE} and \eqref{eq:sin}, respectively. Let $H(q)=\max\left\{ 1+q,\tfrac{3}{2}q \right\}$ and let $\bar{p} \geq 2+4H(q)$ be a sufficiently large even number.
      Then there exist $\gamma > 0$ and $K > 0$ independent of $h$ such that
      \begin{equation*}
            \sup_{ 0 \leq n \leq N }
            \big( \E\big[ |X(t_{n}) - Y_{n}|^{2} \big] \big)^{\frac{1}{2}}
            \leq
            K\big( 1+\big(\E[|X_{0}|^{\bar{p}}]\big)^{\gamma} \big) h.
      \end{equation*}
\end{theorem}

\begin{proof}
Applying techniques in \eqref{eq:tame.emexpect.additive} and \eqref{eq:sine.ed.sub}, we can easily get
\begin{align}
      \label{eq:sine.emexpect.additive}
      \big| \E[X_{t,x}(t+h)-Y^E_{t,x}(t+h)] \big|
      \leq&
      Kh^2(1+|x|^{2+q}),
      \\
      \big| \E[Y^E_{t,x}(t+h)-Y_{t,x}(t+h)] \big|
      \leq&
      Kh^{2}\big(1+|x|^{2+q}\big).\label{eq:sine.ed.additive}
\end{align}
Similarly to \eqref{EulerMean-SquareE.additive} and \eqref{DeviationA.sub}, one can derive that
\begin{align}
      \label{eq:sine.EulerMean-SquareEA}
      \E\big[ |X_{t,x}(t+h)-Y^{E}_{t,x}(t+h)|^{2} \big]
      \leq&
      Kh^{3}\big(1+|x|^{2+q}\big),
      \\
      \label{DeviationA.additive}
      \E\big[|Y^E_{t,x}(t+h)-Y_{t,x}(t+h)|^2\big]
      \leq&
      Kh^{3}(1+|x|^{4+2q}).
\end{align}
Combining \eqref{eq:sine.emexpect.additive}--\eqref{eq:sine.ed.additive} and \eqref{eq:sine.EulerMean-SquareEA}--\eqref{DeviationA.additive} ensures that \eqref{eq:sin} satisfies \eqref{eq:ed}, \eqref{eq:msd} with $p_1=2,p_2=\frac{3}{2}$, which finally completes the proof by Theorem \ref{theorem:convergence}.
\end{proof}

\section{Numerical tests}\label{sec:experiments}

To numerically illustrate the previous theoretical findings, we consider a jump extended version of the $\tfrac32$-volatility model from \cite{Beyn16,Sabanis16}
\begin{equation}\label{eq:numerexa}
\begin{split}
      \mbox{d}X(t)
      =&
      {\mu} X(t^{-})(\nu-|X(t^{-})|)\,\mbox{d}t
      +
      {\xi}|X(t^{-})|^{\frac{3}{2}}\,\mbox{d}W(t)
      \\&+
      {\eta}X(t^{-})\ln(1+X^{2}(t^{-}))\,\mbox{d}\bar{N}(t),
      \quad \forall t \in (0,1], \quad X(0) = 10
\end{split}
\end{equation}
with $\mu = 3 , \nu = 1 , \xi = 0.5, \eta = 0.1$ and an additive noise driven jump-diffusion SDE
\begin{equation}\label{eq:numerexa.additive}
      \mbox{d}X(t)
      =
      (X(t^{-}) - X^{3}(t^{-}))\,\mbox{d}t
      +
      \mbox{d}W(t)
      +
      \mbox{d}\bar{N}(t),
      \quad \forall t \in (0,1],
      \quad X(0)=5.
\end{equation}
Here $\{\bar{N}(t)\}_{0 \leq t \leq 1}$ is a compensated Poisson process with jump intensity $\lambda = 1$. Note that \eqref{eq:numerexa} satisfies Assumptions \ref{as:assumption1}, \ref{as:assumption2} and \ref{ass:drift-polynomial} with polynomial growth rate $q = 2$ and that \eqref{eq:numerexa.additive} obeys Assumptions \ref{as:assumption1}, \ref{as:assumption2} and \ref{as:assumption3} with polynomial growth rate $q = 4$, see Appendix A for more details.
%


To detect the mean-square convergence rates, numerical approximations generated by the tamed (sine) Euler method with a fine stepsize $h=2^{-13}$ are used as the ``exact'' solutions for the order plots.
Then other numerical approximations are calculated by \eqref{eq:tame.method} and \eqref{eq:sin} applied to \eqref{eq:numerexa} and \eqref{eq:numerexa.additive}, respectively,
with five different stepsizes $h = 2^{-i}, i= 8, 9, 10, 11, 12$.
Here the expectations are approximated by the Monte Carlo approximation with $5000$ Brownian and Poisson paths.

\begin{center}
\begin{figure}[!htbp]
    \includegraphics[width=0.475\textwidth]{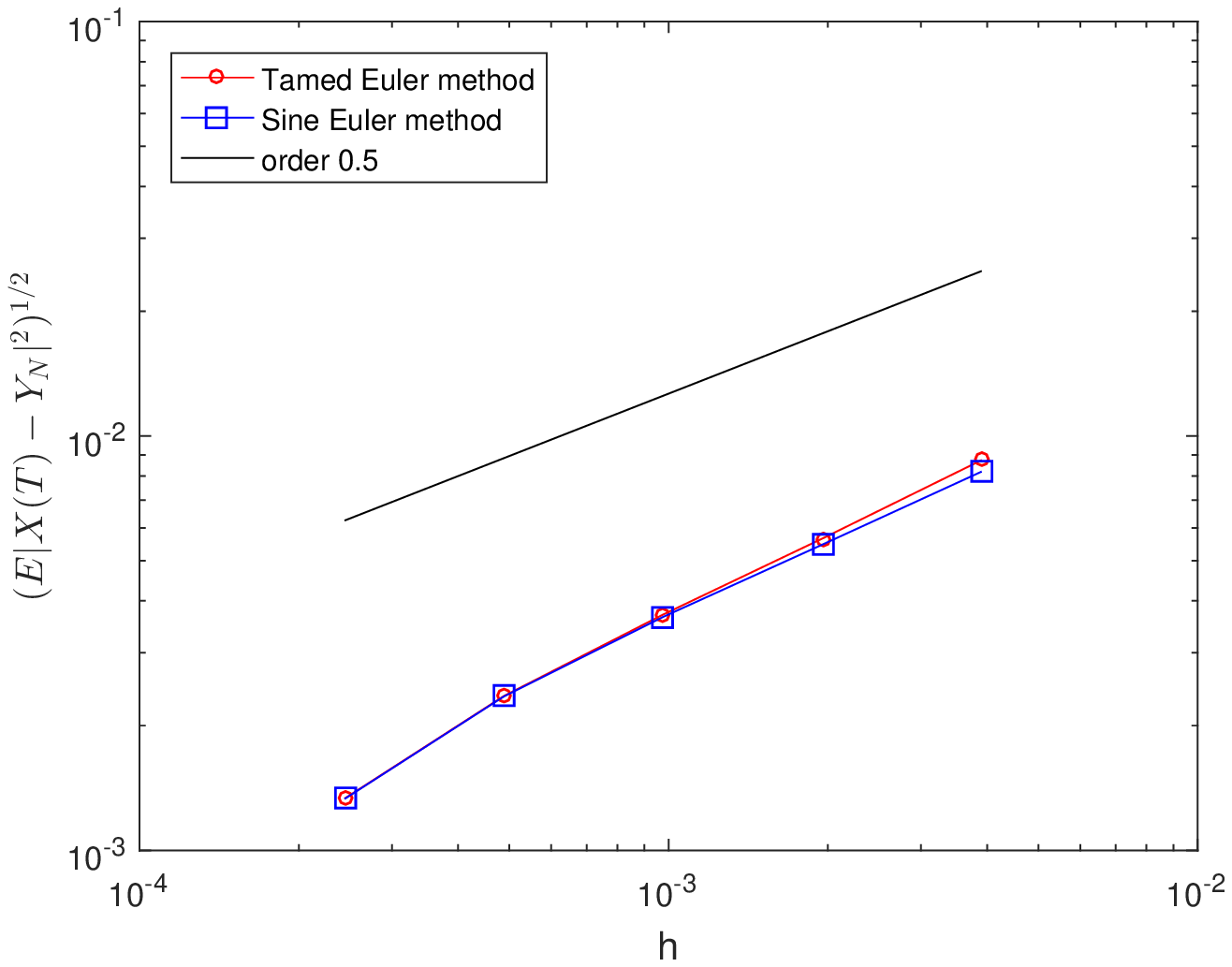}
    \includegraphics[width=0.475\textwidth]{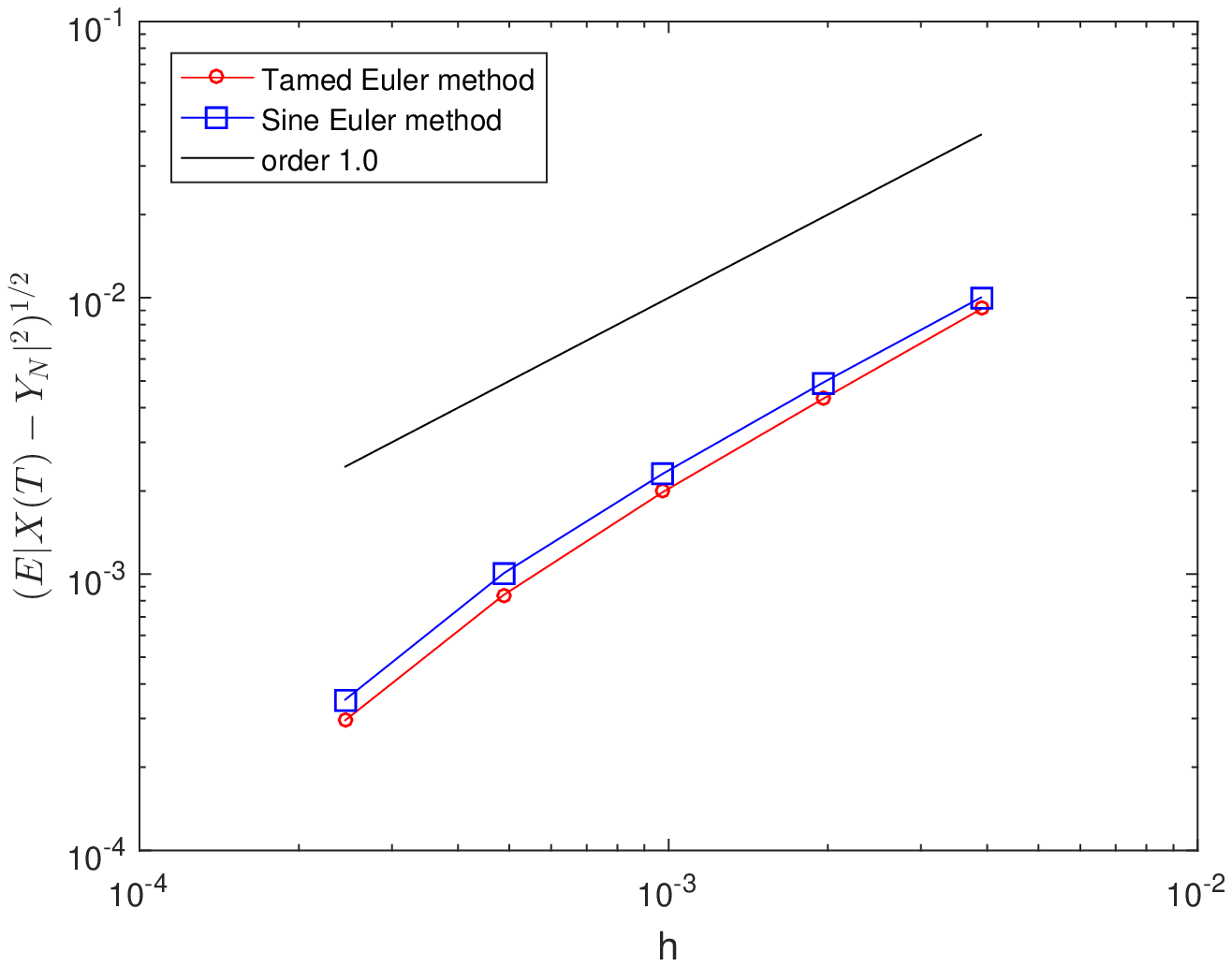}
    \caption{Mean-square convergence rates for \eqref{eq:numerexa} (left) and \eqref{eq:numerexa.additive} (right)
    }
    \label{TSconvergence}
\end{figure}
\end{center}

Figure \ref{TSconvergence} shows that
the slopes of the error lines and the reference lines match well, indicating that the proposed schemes have strong rates of order one-half in non-additive case and order one in additive case.
Additionally, Table 1 lists the CPU time of numerical approximations by \eqref{eq:tame.method} and \eqref{eq:sin}, generated by Matlab R2016a on a desktop (3.86 GB RAM, Intel(R) Core(TM) i5 CPU M480 at 2.67 GHz) with 64 bit Windows 7 operating system.
It seems that the sine Euler method costs slightly less time than the tamed Euler method.

\newcommand{\bb}[1]{\raisebox{-3.5ex}[0pt][0pt]{\shortstack{#1}}}
\begin{center}
\begin{tabular}{c| c c| c c}
  \multicolumn{5}{c}
  {\parbox[c][9mm]{0pt}{} Table 1: CPU time of the tamed and sine methods with different stepsizes
  }\\\hline
  \bb{$h$}
  &
  \multicolumn{4}{c}{\parbox[c][6mm]{0pt}{} CPU time (second)}
  \\\cline{2-5}
  & \multicolumn{2}{c|}{\parbox[c][6mm]{0pt}{} non-additive case}
  & \multicolumn{2}{c}{\parbox[c][6mm]{0pt}{} additive case}
  \\\cline{2-5}
  & \multicolumn{1}{c}{tamed method}
  & \multicolumn{1}{c|}{sine method}
  & \multicolumn{1}{c}{tamed method}
  & \multicolumn{1}{c}{sine method}
  \\\hline
  $2^{-8}$  &0.869440 &0.744448 &0.794577 &0.632519   \\
  $2^{-9}$  &1.203300 &0.932720 &1.031226 &0.776562   \\
  $2^{-10}$ &1.625105 &1.100245 &1.276387 &0.906966   \\
  $2^{-11}$ &2.951017 &1.887072 &2.355608 &1.433305   \\
  $2^{-12}$ &5.789335 &3.588197 &4.325145 &2.473830   \\
  \hline
\end{tabular}
\end{center}

\section*{Appendix A. Verification of assumptions for SDE examples.}
\vspace{0.2em}
\appendix
\renewcommand{\theequation}{A.\arabic{equation}}
\setcounter{equation}{0}

In view of \eqref{eq:numerexa},
the functions $f,g,\sigma \colon \R \to \R$ defined by $f(x) = {\mu} x(\nu-|x|)$, $g(x) = {\xi}|x|^{3/2}$ and $\sigma(x) = {\eta}x\ln(1+x^{2})$ are continuous for all $x \in \R$. Then their derivatives are given by $f'(x) = {\mu}\nu-2{\mu}|x|$, $g'(x) = \frac{3}{2}{\xi}\text{sgn}(x)|x|^{1/2}$ and $\sigma'(x) = \eta(\tfrac{2x^{2}}{1+x^{2}}+\ln(1+x^{2}))$ for all $x \in \R$.
The Appendix in \cite{Sabanis16} tells that
\begin{equation*}
      |f(x)-f(y)|^{2}
      \leq
      3\mu^{2}\max\{1,\nu\}^{2}(1+|x|^{2}+|y|^{2})|x-y|^{2},
      \quad \forall x,y \in \R,
\end{equation*}
which implies that Assumption \ref{ass:drift-polynomial} is satisfied with $q = 2$. This together with Theorems \ref{Theorem:tame.order.I} and \ref{thm:sinerror.sub.nonadditive}
indicates that, for example, $\bar{p} = 20 \geq \max\{2+4H(q),2+2q\}$ is enough for our setting. To verify Assumption \ref{as:assumption1}, we first use the mean value theorem and the H\"{o}lder inequality to get
\begin{equation}\label{eq:testAss31-1}
\begin{split}
      &2\langle x-y,f(x)-f(y) \rangle
      +
      |g(x)-g(y)|^{2}
      +
      \int_{Z}|\sigma(x,z)-\sigma(y,z)|^{2}\,\nu(\mbox{d}z)
      \\\leq&
      |x-y|^{2}
      \int_{0}^{1}
          \big( 2f'(u) + |g'(u)|^{2} + \lambda|\sigma'(u)|^{2} \big)
      \,\mbox{d}r,
      \quad \forall x,y \in \R,
\end{split}
\end{equation}
where $u := y+r(x-y)$. Then the inequality $\ln(1+x^2) \leq 2|x|^{\frac{1}{2}}$ for all $x \in \R$ and \eqref{eq:elementaryinequality} enable us to show
\begin{equation}\label{eq:testAss31-2}
\begin{split}
       &2f'(u) + |g'(u)|^{2} + \lambda|\sigma'(u)|^{2}
       \\\leq&
       2{\mu}\nu
       -
       \big(4{\mu} - \tfrac{9}{4}{\xi}^{2}\big)|u|
       +
       2\lambda\eta^{2}\big( \tfrac{2u^{2}}{1+u^{2}} \big)^{2}
       +
       2\lambda\eta^{2}\big( \ln(1+u^{2}) \big)^{2}
       \\\leq&
       2{\mu}\nu
       +
       8\lambda\eta^{2}
       -
       \big( 4{\mu} - \tfrac{9}{4}{\xi}^{2} - 8\lambda\eta^{2} \big)|u|
       \leq
       2{\mu}\nu
       +
       8\lambda\eta^{2},
\end{split}
\end{equation}
on the condition $4{\mu} - \tfrac{9}{4}{\xi}^{2} - 8\lambda\eta^{2} > 0$.
Hence \eqref{eq:testAss31-1} and \eqref{eq:testAss31-2} prove \eqref{eq:gmc} in Assumption \ref{as:assumption1} for $\mu = 3 , \nu = 1 , \xi = 0.5, \eta = 0.1, \lambda = 1$.  Similarly, we can show \eqref{eq:gcc.a} as follows
\begin{equation*}
\begin{split}
      &2\langle x,f(x) \rangle
      +
      |g(x)|^{2}
      +
      \int_{Z}|\sigma(x,z)|^2\,\nu(\mbox{d}z)
      \\=&
      2{\mu}\nu|x|^{2}
      -
      \big( 2{\mu}-{\xi}^{2} \big)|x|^{3}
      +
      \lambda{\eta}^{2}|x|^{2}(\ln(1+x^{2}))^{2}
      \\\leq&
      2{\mu}\nu|x|^{2}
      -
      \big( 2{\mu}-{\xi}^{2} - 4\lambda{\eta}^{2} \big)|x|^{3}
      \leq
      2{\mu}\nu(1+|x|^2),
      \quad
      \forall x \in \R^d,
\end{split}
\end{equation*}
as $2{\mu}-{\xi}^{2} - 4\lambda{\eta}^{2} > 0$. It remains to verify Assumption \ref{as:assumption2}.
Actually, recalling $\bar{p} = 20$ and using the inequality $\ln(1+x^2) \leq 200 + |x|^{\frac{1}{\bar{p}}}$ for all $x \in \R$ and \eqref{eq:elementaryinequality}, we obtain, after taking $\varepsilon = 1$,
\begin{equation*}
\begin{split}
      &\bar{p}|x|^{\bar{p}-2}\langle x,f(x) \rangle
      +
      \tfrac{\bar{p}(\bar{p}-1)}{2}|x|^{\bar{p}-2}|g(x)|^{2}
      +
      \big(1+(\bar{p}-2) \big)
      \int_Z|\sigma(x,z)|^{\bar{p}}\,\nu(\mbox{d}z)
      \\=&
      \bar{p}{\mu}\nu|x|^{\bar{p}}
      -
      \big( \bar{p}{\mu} - \tfrac{\bar{p}(\bar{p}-1)}{2}{\xi}^{2} \big)|x|^{\bar{p}+1}
      +
      \lambda\big( \bar{p}-1 \big){\eta}^{\bar{p}}
      |x|^{\bar{p}}(\ln(1+x^{2}))^{\bar{p}}
      \\\leq&
      \big(
           \bar{p}{\mu}\nu
           +
           \lambda ( \bar{p}- 1 )(400{\eta})^{\bar{p}}
      \big)|x|^{\bar{p}}
      -
      \big(
           \bar{p}{\mu}
           -
           \tfrac{\bar{p}(\bar{p}-1)}{2}{\xi}^{2}
           -
           \lambda ( \bar{p} - 1 ){(2\eta)}^{\bar{p}}
      \big)
      |x|^{\bar{p}+1}
      \\\leq&
      \bar{p}{\mu}\nu(1+|x|^{\bar{p}}),
      \quad \forall x \in \R^d,
\end{split}
\end{equation*}
as $\bar{p}{\mu}\nu
       +
       \lambda (\bar{p} - 1) (400{\eta})^{\bar{p}}
       > 0$ and
$\bar{p}{\mu} - \tfrac{\bar{p}(\bar{p}-1)}{2}{\xi}^{2}
 - \lambda ( \bar{p}-1  ){(2\eta)}^{\bar{p}} > 0$.
Similarly, we can show that \eqref{eq:numerexa.additive} fulfills Assumptions \ref{as:assumption1}, \ref{as:assumption2} and \ref{as:assumption3} with
$q = 4$.

\section*{Acknowledgments}
The authors are grateful to three anonymous referees whose insightful comments and valuable suggestions are crucial to the improvements of this paper. This paper is dedicated to Prof. Dr. Peter Kloeden in the occasion of his 70th birthday. The third author XW wants to express his gratitude to  Peter for his constant help and encouragement since XW visited the University of Frankfurt am Main, as a joint PhD student.



\providecommand{\href}[2]{#2}
\providecommand{\arxiv}[1]{\href{http://arxiv.org/abs/#1}{arXiv:#1}}
\providecommand{\url}[1]{\texttt{#1}}
\providecommand{\urlprefix}{URL }

\medskip
Received November 2017; 1st revision November 2018; 2nd revision March 2019.\medskip

\end{document}